\documentclass[12pt]{article}
\usepackage{amssymb}
\usepackage{subeqn}
\usepackage{graphicx}
\usepackage{ifpdf}
\ifpdf  \fi
\usepackage{pstricks}
\usepackage{amsfonts,amssymb,slashbox}
\usepackage{mathptmx,helvet,courier,makeidx,multicol,footmisc}
\usepackage[numbers]{natbib}
\bibpunct{(}{)}{;}{a}{,}{,}

\newcommand{\commentout}[1]{}

\newcommand{\ba}{\begin{array}}
        \newcommand{\ea}{\end{array}}
\newcommand{\bc}{\begin{center}}
        \newcommand{\ec}{\end{center}}
\newcommand{\bdm}{\begin{displaymath}}
        \newcommand{\edm}{\end{displaymath}}
\newcommand{\bds} {\begin{description}}
        \newcommand{\eds} {\end{description}}
\newcommand{\ben}{\begin{enumerate}}
        \newcommand{\een}{\end{enumerate}}
\newcommand{\beq}{\begin{equation}}
        \newcommand{\eeq}{\end{equation}}
\newcommand{\bfg} {\begin{figure}[h]}
        \newcommand{\efg} {\end{figure}}
\newcommand{\bi} {\begin {itemize}}
        \newcommand{\ei} {\end {itemize}}
\newcommand{\bqn}{\begin{eqnarray}}
        \newcommand{\eqn}{\end{eqnarray}}
\newcommand{\bqs}{\begin{eqnarray*}}
        \newcommand{\eqs}{\end{eqnarray*}}
\newcommand{\bsl} {\begin{slide}[8.8in,6.7in]}
        \newcommand{\esl} {\end{slide}}
\newcommand{\bss} {\begin{slide*}[9.3in,6.7in]}
        \newcommand{\ess} {\end{slide*}}
\newcommand{\bsq}{\begin{subequations}}
        \newcommand{\esq}{\end{subequations}}
\newcommand{\btb} {\begin {table}}
        \newcommand{\etb} {\end {table}}
\newcommand{\bvb}{\begin{verbatim}}
        \newcommand{\evb}{\end{verbatim}}

\newcommand{\m}{\mbox}
\newcommand {\der}[2] {{\frac {\m {d} {#1}} {\m{d} {#2}}}}

\newcommand {\pd}[2] {{\frac {\partial {#1}} {\partial {#2}}}}

\newcommand{\cas}[1]{{{\left \{ \ba #1 \ea \right. }}}

\newcommand{\reff}[1] {{{Figure \ref {#1}}}}
\newcommand{\refe}[1] {{(\ref {#1})}}

\def\la      {{\lambda}}

\def\pmb#1{\setbox0=\hbox{$#1$}%
   \kern-.025em\copy0\kern-\wd0
   \kern.05em\copy0\kern-\wd0
   \kern-.025em\raise.0433em\box0 }

\def\r{{\rho}}

\def\dx     {{\Delta x}}
\def\dt     {{\Delta t}}

\begin{document}
\title{Control of a lane-drop bottleneck through variable speed limits} 
\author{Hui-Yu Jin \footnote{Department of Automation, Xiamen University, Xiamen 361005, P.~R.~China. Email: jinhy@xmu.edu.cn} ~ and Wen-Long Jin \footnote{Department of Civil and Environmental Engineering, California Institute for Telecommunications and Information Technology, Institute of Transportation Studies, 4000 Anteater Instruction and Research Bldg, University of California, Irvine, CA 92697-3600. Tel: 949-824-1672. Fax: 949-824-8385. Email: wjin@uci.edu. Corresponding author}}
\maketitle

\begin{abstract}
The discharging flow-rate of a lane-drop bottleneck can drop when its upstream is congested, and such capacity drop can lead to additional traffic congestion as well as safety threats. Even though many studies have demonstrated that variable speed limits (VSL) can effectively delay and even avoid the occurrence of capacity drop, there lacks a simple approach for analyzing the performance of a control system. In this study, we formulate the VSL control problem for the traffic system in a zone upstream to a lane-drop bottleneck based on two traffic flow models: the Lighthill-Whitham-Richards (LWR) model, which is an infinite-dimensional partial differential equation, and the link queue model, which is a finite-dimensional ordinary differential equation. In both models, the discharging flow-rate is determined by a recently developed model of capacity drop, and the upstream in-flux is regulated by the speed limit in the VSL zone. Since the link queue model approximates the LWR model and is much simpler, we first analyze the control problem and develop effective VSL strategies based on the former. First for an open-loop control system with a constant speed limit, we prove that a constant speed limit can introduce an uncongested equilibrium state, in addition to a congested one with capacity drop, but the congested equilibrium state is always exponentially stable. Then we apply a feedback proportional-integral (PI) controller to form a closed-loop control system, in which the congested equilibrium state and, therefore, capacity drop can be removed by the I-controller. Both analytical and numerical results show that, with appropriately chosen controller parameters, the closed-loop control system is stable, effect, and robust. Finally,  we show that the VSL strategies based on I- and PI-controllers are also stable, effective, and robust for the LWR model. Since the properties of the control system are transferable between the two models, we establish a dual approach for studying the control problems of nonlinear traffic flow systems. We also confirm that the VSL strategy is effective only if capacity drop occurs. The obtained method and insights can be useful for future studies on other traffic control methods and implementations of VSL strategies in lane-drop bottlenecks, work zones, and incident areas.          
\end{abstract}


\section{Introduction}

An important characteristic of various freeway bottlenecks is that, when a bottleneck becomes active; i.e., when one of the upstream branches is congested but the downstream branch is not, the maximum discharging flow-rate of the bottleneck can drop to a lower level than that in free flow.
In other words, there can be two distinctive capacities for uncongested and congested traffic at such bottlenecks. This is the so-called two-capacity or capacity-drop phenomenon of active bottlenecks, in which ``maximum flow rates decreases when queues form''  \citep{banks1990flow,banks1991twocapacity,hall1991capacity}.
Such capacity drop has been observed at merges, tunnels, lane drops, curves, and upgrades, where the bottlenecks cannot provide sufficient space for upstream vehicles \citep{chung2007relation,kim2012capacity}. In free flow, the bottleneck's capacity can be fully utilized, and the maximum flow-rate for the bottlenecks can reach around 2300 vphpl, which is about the same as the capacity of a homogeneous road lane without bottlenecks \citep{hcm1985,hall1991capacity}. However, with excessive traffic demands the traffic system can break down, the maximum discharge flow-rate drops below the free-flow capacity. It was shown that such dropped capacity is stable, although interactions among several bottlenecks can cause fluctuations in discharging flow-rates \citep{kim2012capacity}.
Depending on the geometry of a bottleneck, the magnitude of capacity drop is in the order of 10\% \citep{persaud1998exploration,cassidy1999bottlenecks, bertini2005empirical,chung2007relation}.
In \citep{persaud1998exploration,persaud2001breakdown}, traffic breakdown and capacity drop were found to be related to the upstream traffic demand randomly.
Also closely associated with capacity drop is the existence of discontinuous flow-density relations in fundamental diagrams. For examples, in \citep{edie1961cf}, discontinuity was observed in the flow-density relation for the Lincoln tunnel; in \citep{drake1967fd,koshi1983fd}, a flow-density relation of the reverse-lambda shape was calibrated; more studies have confirmed the existence of discontinuous fundamental diagrams \citep{payne1984discontinuity,hall1992fd}.
In \citep{hall1991capacity}, it was shown that such discontinuous fundamental diagrams generally arise in the congested area of an active bottleneck and suggested that the discontinuity is caused by the capacity drop phenomenon.

In the literature, capacity drop is modeled by discontinuous fundamental diagrams \citep{lu2008explicit,lu2009entropy}. However, empirical observations suggest that fundamental diagrams are still continuous in steady-state traffic flow \citep{cassidy1998bivariate}, and theoretical analyses show that unrealistic, infinite  characteristic wave speeds can arise from a discontinuous flow-density relation   \citep{li2013modeling}. In contrast, in \citep{jin2013_cd} a new kinematic wave model of capacity drop was proposed. In the new model, fundamental diagrams are still continuous, but the entropy condition for uniquely solving the system of hyperbolic conservation laws is discontinuous. The new entropy condition is a junction flux function used to prescribe discharging flow-rates from the upstream demand and downstream supply. Therefore, the model can be easily discretized and incorporated into the Cell Transmission Model (CTM) \citep{daganzo1995ctm,lebacque1996godunov}. Furthermore, the model was shown to be well defined analytically, and the new junction flux function was shown to be invariant \citep{lebacque2005network,jin2012_riemann}. Therefore, the junction flux can also be incorporated into the link queue model \citep{jin2012_link}, which is a system of ordinary differential equations but approximates the kinematic wave model. In addition, with empirical observations it was shown that the new model can explain the existence of discontinuous flow-density relation in steady traffic states.  

Clearly, a drop at the active bottleneck's capacity can reduce the total discharging flow-rate of the whole corridor and prolong vehicles' travel times \citep{daganzo1999remarks}.
That the capacity of a road network may drop substantially when it is most needed during the peak period has been a baffling nature of freeway traffic dynamics \citep{papageorgiou2002freeway}.
Therefore, an important motivation and theoretical foundation for developing ramp metering, variable speed limits, and other control strategies is to avoid or delay the occurrence of capacity drop \citep{banks1991metering,papageorgiou1991alinea,papageorgiou1997ALINEA,cassidy2005merge,papageorgiou2005review,papageorgiou2007alinea}.
In \citep{papageorgiou2008effects}, it was shown that mainline freeway¡¯s variable speed limits (VSL) can impact fundamental diagrams and could increase discharging flow-rates in congested traffic. In \citep{allaby2006variable}, microscopic simulations based on PARAMICS are used to study the impacts of VSL on both safety and efficiency of the eastbound Queen Elizabeth Way (QEW) located near Toronto, Canada, but it showed that travel times can be increased by VSL, due to the temporal and spatial limitations of VSL algorithms. In \citep{cascetta2011empirical}, it was shown that speed enforcement can reduce bottleneck congestion and travel time with empirical evidence. In \citep{hegyi2005model}, network-wide coordination of VSL and ramp metering was studied within the framework of model predictive control based on a higher-order simulation model.  In \citep{lu2010new}, an integrated variable speed limit and ramp metering strategy was proposed, and it was shown that by creating a discharging section right before a bottleneck can increase discharging flow-rate. In \citep{carlson2010optimal}, with a macroscopic second-order traffic flow model, it was shown that VSL have similar effects as ramp metering; coordination of VSL control and ramp metering is solved as a  discrete-time optimal control problem for large-scale networks.

In this study, we study the variable speed limit (VSL) strategy to tackle the capacity drop at a lane-drop bottleneck as shown in \reff{lanedrop_vsl_ill}, which can be caused by lane closures, work zones, or incidents. 
The capacity in the downstream section, $C$, is determined by systematic lane changes inside zones 1 and 2 \citep{jin2010lc,jin2010lcchar,jin2012_lc,gan2013validation}.
When the traffic demand in zone 0, $d(t)$, is higher than $C$, a queue forms in zone 0, the capacity drops by $\Delta$ to $C(1-\Delta)$, and vehicles accelerate to the free-flow speed in zones 1 and 2. 
However, when the VSL zone is introduced upstream to zone 0, its demand, $d^-(t)$, can be controlled by the variable speed limit, $u(t)$. Then $d(t)$ in zone 0 can be restricted, so that zones 1 and 2 will be clear and the total discharging flow-rate can reach $C$. 

\bfg\bc
\includegraphics[width=6in]{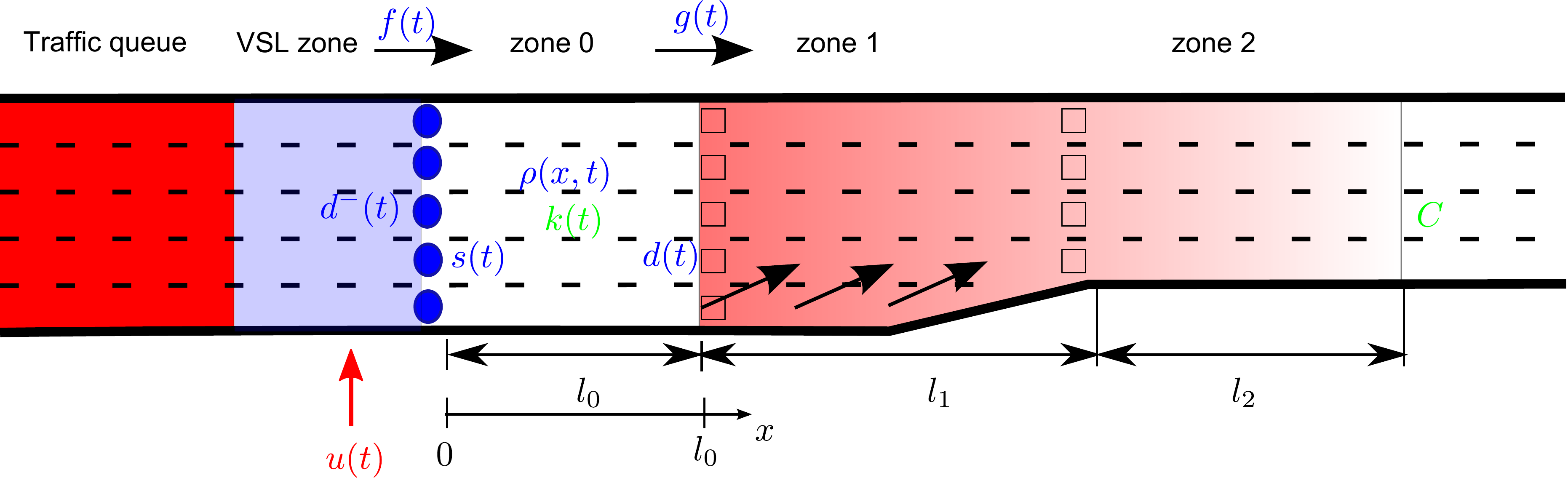}\caption{An illustration of variable speed limits at a lane-drop bottleneck}\label{lanedrop_vsl_ill}
\ec\efg

As the traffic queue in zone 0 activates capacity drop, we will study the VSL strategy as a feedback control problem of traffic dynamics in zone 0. For zone 0, the traffic dynamics can be described by the Lighthill-Whitham-Richards (LWR) model \citep{lighthill1955lwr,richards1956lwr}, its in-flux calculated from $d^-(t)$ and $s(t)$ with a traditional flux function consistent with the Lax entropy condition \citep{jin2009sd}, and its out-flux calculated from $d(t)$ and $C$ with the new capacity drop model in \citep{jin2013_cd}.
Therefore, the goal of this study is to solve $u(t)$ from the traffic density near the boundary between zones 0 and 1. 
Instead of directly solving the control problem with the LWR model, which is a partial differential equation, we first solve the problem with the link queue model, which is an ordinary differential equation.

The rest of the paper is organized as follows. In Section 2, we present two models for traffic dynamics in zone 0: the LWR model and the link queue model, for which the speed limit restricts the upstream in-flux and the capacity drop determines the discharging flow-rate (downstream out-flux). In Section 3, with the link queue model, we analyze the equilibrium states and their stability of the open-loop control system with different constant speed limits. In Section 4, we consider a closed-loop control system with an I- or PI-controller and analyze the equilibrium states and their stability. In Section 5, with numerical simulations, we verify the stability of the closed-loop system and examine the robustness of the designed VSL controller subject to disturbances in the target density and the demand pattern. In Section 6, with numerical simulations, we show that both I- and PI-controllers are stable, effective, and robust for the LWR model. Finally in Section 7 we conclude this study. 

\section{Two models of the traffic system}
In this study, the traffic system inside zone 0 in \reff{lanedrop_vsl_ill} is considered as the control object, and traffic dynamics inside other zones are not explicitly modeled.
In this section,  we present two models of the traffic system: the LWR model \citep{lighthill1955lwr,richards1956lwr} and the link queue model \citep{jin2012_link}. In particular, the speed limit in the VSL zone, $u(t)$, is introduced to regulate the upstream in-flux, and a novel capacity drop model in \citep{jin2013_cd} is used to determine the discharging flow-rate.

We assume that $x\in [0, l_0]$ and $x$ increases in the traffic direction. We denote the traffic demand to the whole system at $x=0^-$ by $d^-(t)$. We denote upstream supply and downstream demand of zone 0 by $s(t)$ and $d(t)$, respectively. Here the demand and supply concepts were first introduced in \citep{lebacque1996godunov} and are the same as sending and receiving flows in \citep{daganzo1995ctm}. In addition, we denote the in- and out-fluxes by $f(t)$ and $g(t)$, respectively.

The state variable of the traffic system inside zone 0 depends on the used traffic flow model and will be introduced in the following subsection. Conceptually, the system dynamics are related to the initial conditions and boundary conditions at $x=0$ and $x=l_0$; the control signal, $u(t)$, $d^-(t)$, and $s(t)$ will regulate the in-flux; and the capacity drop model determines the out-flux, i.e., discharging flow-rate, $g(t)$. Then to solve the control problem we need to design an algorithm to calculate $u(t)$ from observed traffic conditions inside zone 0.

\subsection{The LWR model}
In the LWR model, the state variable is the traffic density $\r(x,t)$ at location $x$ and time $t$. Related variables are speed $v(x,t)$ and flow-rate $q(x,t)$. From the conservation of traffic flow, $\pd{\r}t+\pd{q}x=0$, and the fundamental diagram, $q=Q(\r)$, which is unimodal, we have the following evolution equation \citep{lighthill1955lwr,richards1956lwr}:
\bqn
\pd{\r} t+\pd{Q(\r)}x&=&0, \label{lwr}
\eqn
which is a nonlinear hyperbolic conservation law \citep{smoller1983shock}. Note that $q=Q(\r)$ is continuous even when capacity drop occurs \citep{cassidy1998bivariate,jin2013_cd}. Many empirical observations and car-following models suggest that the following triangular fundamental diagram is an appropriate approximation of the flow-density relation in steady states \citep{munjal1971multilane,haberman1977model,newell1993sim}:
\bsq \label{tri-fd}
\bqn
Q(\r)=\min\{v_f \r, w(k_j-\r)\},
\eqn
where $v_f$ is the free-flow speed, $w$ the shock wave speed in congested traffic, $k_j$ the jam density for all lanes. Therefore the maximum flow-rate, i.e., the capacity, is attained at the critical density:
\bqn
k_c&=&\frac{w}{v_f+w}k_j.
\eqn
\esq

In \citep{jin2009sd}, it was shown that an entropy condition for uniquely solving the LWR model, \refe{lwr}, is given by 
\bqn
q(x,t)&=&\min\{\delta(x^-,t),\sigma(x^+,t)\}, \label{entropycond}
\eqn
where the traffic demand and supply are defined respectively by ($0<x<l_0$)
\bsq \label{d-s}
\bqn
\delta(x,t)&=&D(\r(x,t))\equiv Q(\min\{\r(x,t),k_c\})=\min\{v_f k_c, v_f \r(x,t)\},\\
\sigma(x,t)&=&S(\r(x,t))\equiv Q(\max\{\r(x,t),k_c\})=\min\{v_f k_c, w(k_j-\r(x,t))\}.
\eqn
\esq
Note that, in the continuous flux function, \refe{entropycond},  $\delta(0^-,t)$ and $\sigma(l_0^+,t)$ cannot be calculated from the state variables. Thus $f(t)\equiv q(0,t)$ and $g(t)\equiv q(l_0,t)$ have to be determined by boundary conditions at $x=0$ and $x=l_0$ in the LWR model.

At the downstream boundary $x=l_0$, we apply the phenomenological capacity drop model \citep{jin2013_cd}:
\bqs
g(t)&=&\min\{d(t), C(1-\Delta \cdot I_{d(t)>C})\},
\eqs
where $d(t)=\delta(l_0^-,t)$, $\Delta$ is the magnitude of capacity drop, and $I_{d(t)>C}$ is the indicator function equal to 1 when $D(k)>C$ and 0 otherwise.
Since $C<v_f k_c$, the discharging flow-rate can be simplified as
\bqn
g(t)&=&\min\{v_f \r(l_0^-,t), C(1-\Delta \cdot I_{\r(l_0^-,t)>k_1})\}=\cas{{ll} v_f \r(l_0^-,t), & \r(l_0^-,t)\leq k_1\\ C(1-\Delta), & \r(l_0^-,t)> k_1} \label{dischargingflow-lwr}
\eqn
where $k_1= C/v_f$. Obviously this discharging flow-rate is discontinuous at $k_1$ and drops from $C$ to $C(1-\Delta)$ when $\r(l_0^-,t)$ exceeds $k_1$. This captures the phenomenon of capacity drop: when the demand at $x=l_0^-$ is higher than the downstream capacity, $C$, a queue forms in zone 0, and the discharging flow-rate drops below the downstream capacity. Note that, in the discrete system, \refe{ctm}, $\r(l_0^-,t)$ is replaced by $\r_n^j$, the density in the last cell.

We assume that the variable speed limit is applied at the upstream boundary ($x=0^-$) and $u(t)\leq v_f$.
Then the in-flux is bounded by $\frac{u(t)}{u(t)+w} wk_j<v_f k_c$. That is, we extende \refe{entropycond} at $x=0$ by 
\bqn
f(t)&=&\min\{d^-(t), \frac{u(t)}{u(t)+w} wk_j, w(k_j-\r(0^+,t))\}, \label{infloweqn-lwr}
\eqn
where $s(t)=\sigma(0^+,t)=\min\{v_f k_c, w(k_j-\r(0^+,t))\}$.
Obviously, if there is no VSL control; i.e., if $u(t)=v_f$, \refe{infloweqn-lwr} is equivalent to \refe{entropycond}. Note that, in the discrete system, \refe{ctm}, $\r(0^+,t)$ is replaced by $\r_1^j$, the density in the first cell.

Therefore, the LWR model, \refe{lwr}, the entropy condition for $0<x<l_0$, \refe{entropycond}, the capacity drop model at $x=l_0$, \refe{dischargingflow-lwr}, and the VSL model at $x=0$, \refe{infloweqn-lwr}, form a control system, which has an infinite number of state variables, $\r(x,t)$, captures the phenomenon of capacity drop at $x=l_0$, and has a control signal $u(t)$. The corresponding discrete system is also formulated with \refe{ctm}.

The triangular fundamental diagram in \refe{tri-fd} and related variables are illustrated in \reff{vsl_fd}, in which $v_1= \frac{C}{k_jw-C} w$, $v_2= \frac{C(1-\Delta)}{k_jw-C(1-\Delta)} w$, and $v_2<v_1<v_f$. A traffic system is congested if $k>k_c$ and uncongested otherwise.

\bfg\bc
\includegraphics[width=5in]{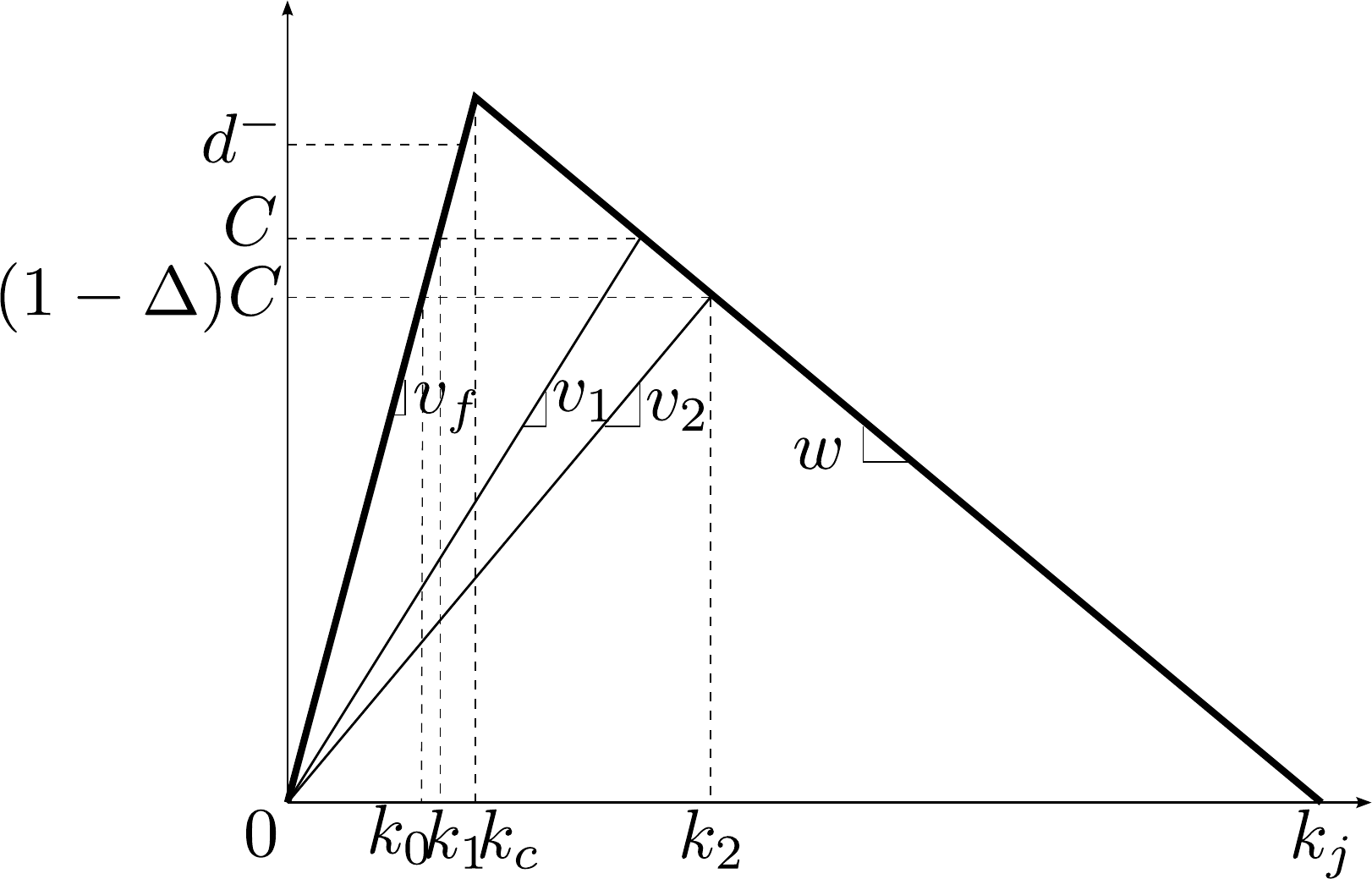}
\caption{The fundamental diagram and variables}\label{vsl_fd}
\ec\efg

\subsection{The link queue model}
The LWR model has been useful for analyzing the formation and dissipation of traffic queues on a road through shock and rarefaction waves. However, it is an infinitely dimensional system, and the control problem of such a partial differential equation is still an open problem. In contrast, there exist many tools for studying the control problem of ordinary differential equations, even though it is still quite challenging to control a nonlinear system.

Fortunately, in \citep{jin2012_link}, it was shown that the LWR model can be approximated by a so-called link queue model, which is a nonlinear ordinary differential equation. The new system describes the dynamics in the average density of zone or link 0, denoted by $k(t)$ here, based on the following conservation equation
\bqn\label{control_1}
\der{k(t)}t&=&\frac 1 {l_0} (f(t) -g(t)  ).
\eqn
In addition, the demand and supply function in \refe{d-s} can be directly extended for the link density:
\bsq
\bqn
d(t)&=&\min\{v_f k, v_f k_c\},\\
s(t)&=&\min\{v_f k_c, w(k_j-k)\},
\eqn
\esq
where $d(t)$ and $s(t)$ are the demand and supply of the whole link.
However, the in- and out-fluxes, $f(t)$ and $g(t)$, have to be determined by the upstream and downstream boundary conditions.

The phenomenological capacity drop model \citep{jin2013_cd} can also be used in the link queue model to determine the out-flux by revising \refe{dischargingflow-lwr} as 
\bqn
g(t)&=&\min\{v_f k(t), C(1-\Delta \cdot I_{k(t)>k_1})\}=\cas{{ll} v_f k(t), & k(t)\leq k_1\\ C(1-\Delta), & k(t)> k_1} \label{dischargingflow}
\eqn
This also captures the phenomenon of capacity drop: when the demand in zone 0 is higher than the downstream capacity, $C$, a queue forms in zone 0, and the discharging flow-rate drops below the downstream capacity.

Similarly, the in-flux can be written as
\bqn
f(t)&=&\min\{d^-(t), \frac{u(t)}{u(t)+w} wk_j, w(k_j-k(t))\}. \label{infloweqn}
\eqn

Since both $f(t)$ and $g(t)$ are functions of $k(t)$ and $u(t)$; i.e, $f(t)=F(u(t),k(t))$, and $g(t)=G(k(t))$, \refe{control_1}, \refe{dischargingflow}, and \refe{infloweqn} form a simple control system, 
\bqn\label{system_8au}\label{controlsys}
\der{}{t} k &=& \frac 1 {l_0} (F(u,k)-G(k))\\
&=&\frac 1 {l_0} (\min\{d^-(t), \frac{u(t)}{u(t)+w} wk_j, w(k_j-k(t))\}-\min\{v_f k(t), C(1-\Delta \cdot I_{k(t)>k_1})\} ).\nonumber
\eqn
in which $k(t)$ is the state variable, $u(t)$ the control variable, and the capacity drop phenomenon is captured.
It is noticeable that the variable $t$ does not appears in  the right side of \refe{system_8au}. Such a system is called an autonomous system in control theory. 

The link queue model in \refe{controlsys} is still highly nonlinear, since the right-hand side is discontinuous. However, compared with the system based on the LWR model, this system is much simpler with only one state variable. Therefore, we will first analytically study the control system and design the control algorithms based on the link queue model and then verify the results with the LWR model.

\section{An open-loop control system with a constant speed limit}
In this section, we assume that the traffic demand $d^-(t)=d^-$ is constant\footnote{Note that when $d^-(t)$ varies with time, the traffic patterns will be more complicated, but the insights can still be applied.}. With a constant speed limit,  we analyze the equilibrium states and their stability for the open-loop control system based on the link queue model, \refe{system_8au}.

\subsection{Equilibrium states}
When the control system reaches an equilibrium state, in which $k(t)=k^*$ and $u(t)=u^*$, we have $\der{}{t}k(t)=0$, and 
\bqn
F(u^*,k^*)=G(k^*). \label{equpoint}
\eqn  
In the equilibrium state, the in- and out-fluxes are the same at $g^*$.

In equilibrium states \refe{equpoint} can be written as
\bqs
\min\{d^-, \frac{u^*}{u^*+w} wk_j, w(k_j-k^*)\}=\min\{v_f k^*, C(1-\Delta \cdot I_{k^*>k_1})\}.
\eqs
As shown in \reff{vsl_equilibrium}, the equilibrium states are the intersection points between the thick solid curve ($G(k)$) and the thick dotted curve ($F(u,k)$). 

\bfg\bc
\includegraphics[width=5in]{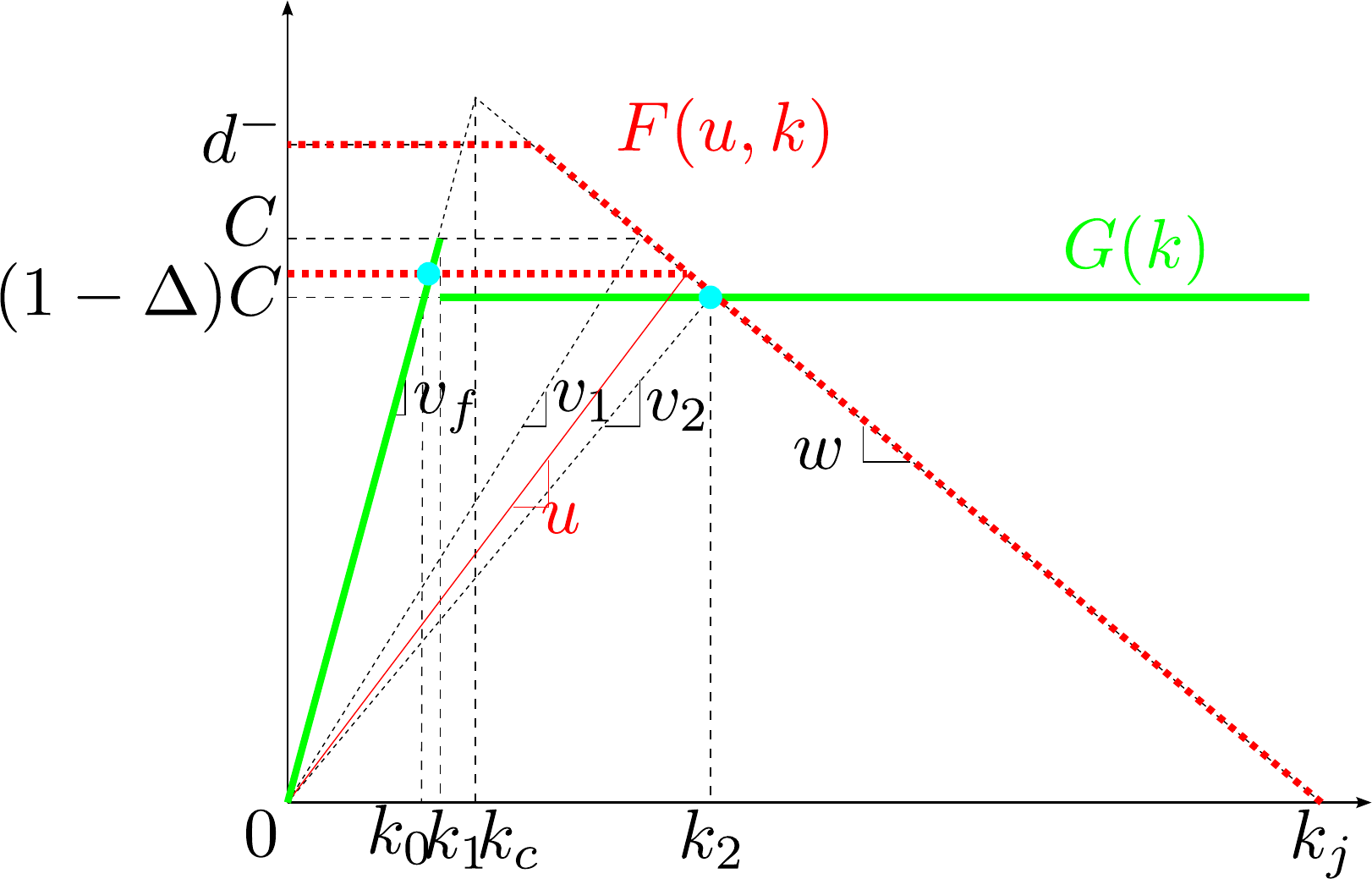}
\caption{Equilibrium states in the fundamental diagram}\label{vsl_equilibrium}
\ec\efg

First, if we set $u^*=v_f$, then no VSL control is applied. In this case, we have the following equilibrium states.
\ben
\item When $d^-> C$, there is only one congested equilibrium state at $k^*=k_2$.
\item When $(1-\Delta) C \leq d^- \leq C$, there are two equilibrium states: uncongested at $k^*=\frac{d^-}{v_f}$ and congested at $k^*=k_2$.
\item When $d^-<(1-\Delta)C$, there is only one uncongested equilibrium state at $k^*=\frac{d^-}{v_f}$.
\een

Second, if we introduce VSL control, i.e., if $u^*<v_f$, we have the following equilibrium states.
\ben
\item When $d^->C$, the traffic system can have multiple equilibrium states, depending on $u^*$ and the initial state $k(0)$.
\ben
\item When $u^*> v_1$,  $\frac{u^*}{u^*+w} wk_j>C$, and  the equilibrium state is always congested at $k^*=k_2$ and $g^*=(1-\Delta)C$ for $k(0)\leq k_1$.
\item When $u^*< v_2$,  $\frac{u^*}{u^*+w} wk_j<(1-\Delta)C$, and the equilibrium state is always uncongested at $k^*=\frac{u^*}{u^*+w} \frac{w}{v_f} k_j<\frac{(1-\Delta)C}{v_f}$ and $g^*=\frac{u^*}{u^*+w} wk_j$ for any $k(0)$. 
\item When $v_2\leq u^*\leq v_1$, $(1-\Delta)C \leq \frac{u^*}{u^*+w} wk_j\leq C$, and the equilibrium state can be uncongested at $k^*=\frac{u^*}{u^*+w} \frac{w}{v_f} k_j$ and $g^*=\frac{u^*}{u^*+w} wk_j$ for $k(0)\leq k_1$, or congested at $k^*=k_2$ and $g^*=(1-\Delta)C$ for $k(0)>k_1$.
\een
Therefore, when $d^->C$, the optimal speed limit is $u^*=v_1$, and the optimal discharging flow-rate is $C$ if $k(0)\leq k_1$.

\item When $d^-\in((1-\Delta)C, C]$, the traffic system can have multiple equilibrium states, depending on $u^*$ and the initial state.
\ben
\item When $u^*\geq \frac{d^-}{k_jw-d^-}w$, $\frac{u^*}{u^*+w} wk_j\geq d^-$, and the equilibrium state can be uncongested at $k^*=\frac{d^-}{v_f}$ and $g^*=d^-$ for $k(0)\leq k_1$, or congested at $k^*=k_2$ and $g^*=(1-\Delta)C$.
\item When $u^*< v_2$,  $\frac{u^*}{u^*+w} wk_j<(1-\Delta)C$, and the equilibrium state is always uncongested at $k^*=\frac{u^*}{u^*+w} \frac{w}{v_f} k_j<\frac{(1-\Delta)C}{v_f}$ and $g^*=\frac{u^*}{u^*+w} wk_j<(1-\Delta)C$ for any $k(0)$.
\item When $v_2\leq u^*< \frac{d^-}{k_jw-d^-}w$, $(1-\Delta)C \leq\frac{u^*}{u^*+w} wk_j< d^-$, and the equilibrium state can be uncongested at $k^*=\frac{u^*}{u^*+w} \frac{w}{v_f} k_j$ and $g^*=\frac{u^*}{u^*+w} wk_j<d^-$ for $k(0)\leq k_1$, or congested at $k^*=k_2$ and $g^*=(1-\Delta)C$ for $k(0)>k_1$.
\een
Therefore, when $d^-\in((1-\Delta)C, C]$, the optimal speed limit is $u^*\geq \frac{d^-}{k_jw-d^-}w$, and the optimal discharging flow-rate is $g^*=d^-$ if $k(0)\leq k_1$.

\item When $d^-= (1-\Delta)C$, the traffic system can have multiple equilibrium states, depending on $u^*$ and the initial state.
\ben 
\item When $u^*\geq \frac{d^-}{k_jw-d^-}w=v_2$, $\frac{u^*}{u^*+w} wk_j\geq d^-$, and the equilibrium state can be uncongested at $k^*=\frac{d^-}{v_f}$ and $g^*=d^-$, if $k(0)\leq k_1$, or congested at $k^*=k_2$ and $g^*=(1-\Delta)C=d^-$, if $k(0)>k_1$.
\item When $u^*< v_2$,  $\frac{u^*}{u^*+w} wk_j<d^-$, and the equilibrium state is always uncongested at $k^*=\frac{u^*}{u^*+w} \frac{w}{v_f} k_j<\frac{(1-\Delta)C}{v_f}$ and $g^*=\frac{u^*}{u^*+w} wk_j<(1-\Delta)C$ for any $k(0)$.
\een
Therefore, when $d^-= (1-\Delta)C$, the optimal speed limit is $u^*\geq v_2$, and the optimal discharging flow-rate is $g^*=d^-$ for any $k(0)$.

\item When $d^-< (1-\Delta)C$, the traffic system can have multiple equilibrium states, depending on $u^*$ .
\ben 
\item When $u^*\geq \frac{d^-}{k_jw-d^-}w$, $\frac{u^*}{u^*+w} wk_j\geq d^-$, and the equilibrium state is always uncongested at $k^*=\frac{d^-}{v_f}$ and $g^*=d^-$ for any $k(0)$.
\item When $u^*< \frac{d^-}{k_jw-d^-}w$,  $\frac{u^*}{u^*+w} wk_j<d^-$, and the equilibrium state is always uncongested at $k^*=\frac{u^*}{u^*+w} \frac{w}{v_f} k_j<\frac{(1-\Delta)C}{v_f}$ and $g^*=\frac{u^*}{u^*+w} wk_j<(1-\Delta)C$ for any $k(0)$.
\een
Therefore, when $d^-< (1-\Delta)C$, the optimal speed limit is $u^*\geq \frac{d^-}{k_jw-d^-}w$, and the optimal discharging flow-rate is $g^*=d^-$ for any $k(0)$.
\een

From the analysis above, we have the following observations:
\ben
\item No VSL control is needed when $d^-\leq (1-\Delta)C$. That is, we can set $u^*=v_f$, which can lead the optimal discharging flow-rate $d^-$ under any initial conditions.
\item However, comparing without VSL control, the constant VSL controller successfully introduces one uncongested equilibrium state when $d^->C$ with the maximum discharging flow-rate. This justifies the introduction of the VSL control. 
\item The optimal speed limit is always $u^*\geq \frac{d^-}{k_jw-d^-}w$, when $d^-\leq C$; $u^*=v_1$ when $d^->C$.
\item The optimal discharging flow-rate is $g^*=\min\{C,d^-\}$. But when $d^->(1-\Delta)C$, the optimal discharging flow-rate is only achieved when $k(0)\leq k_1$.
\een
Therefore, when $d^->(1-\Delta)C$, the VSL control algorithm has to achieve two goals: to pick an optimal speed limit, and to drive the initial condition $k(0)\leq k_1$. 

In addition, if there is no capacity drop; i.e., if $\Delta=0$, then $v_1=v_2$, and the optimal discharging flow-rate $g^*=\min\{C, d^-\}$ is achieved, if and only if $u^*\geq \frac{d^-}{k_jw-d^-}w$. Since $v_f \geq \frac{d^-}{k_jw-d^-}w$, it means that there is no need to implement VSL when there is no capacity drop. This analysis corroborates that the key contribution of VSL is to avoid capacity drop.

\subsection{Stability analysis of the control system with a constant speed limit}
In the preceding subsection, we found that $u^*=v_1$ is an optimal speed limit when $d^->(1-\Delta)C$. That is, we can achieve the optimal discharging flow-rate of $g^*=\min\{C,d^-\}$ with this speed limit. Then can we simply set the speed limit as $v_1$ in the VSL zone? If so, this will lead to very simple implementations of VSL strategy in the real world.

To answer this question, we first note that there are two equilibrium states when $d^->(1-\Delta)C$ and $u^*=v_1$: uncongested at $k^*=\frac{\min\{C,d^-\}}{v_f}$ and $g^*=\min\{C,d^-\}$, or congested at $k^*=k_2$ and $g^*=(1-\Delta)C$. In this subsection, we analyze the stability of these equilibrium states, which will be critical to answer the question above. 

When $d^->C$, we set $u(t)=v_1$ and  rewrite \refe{system_8au} as
\bqn
\der{}{t}k&=&\frac 1{l_0}(\min\{C, w(k_j-k)\}-\min\{v_f k, (1-\Delta \cdot I_{k>k_1})C\}). \label{syssta}
\eqn
\ben
\item For the uncongested equilibrium state at $k^*=k_1$ and $g^*=C$, we denote a small disturbance in traffic density by $z=k-k_1$, then \refe{syssta} is equivalent to
\bqs
\der{}{t}z&=&-\frac 1{l_0}\min\{v_f z, -\Delta \cdot I_{z>0} \cdot C\}=\cas{{ll} -\frac {v_f}{l_0}z, &z\leq 0\\ \frac {C\Delta}{l_0}, & z>0} 
\eqs
Therefore $\der{z}{t}\geq 0$, and $z=0$ is a unstable saddle point. That is, $k^*=k_1$ is unstable if the disturbance is positive.

\item For the congested equilibrium state at $k^*=k_2$ and $g^*=(1-\Delta)C$, we denote a small disturbance in traffic density by $z=k-k_2$, then \refe{syssta} is equivalent to
\bqs
\der{}{t}z&=&-\frac w{l_0} z.
\eqs
Therefore $z=0$ is an exponentially stable equilibrium point, and $k^*=k_2$ is exponentially stable.
\een
Thus when $d^->C$, with $u^*=v_1$, any positive disturbance in density would drive the system away from the optimal equilibrium state, and the discharging flow-rate will be reduced to $(1-\Delta)C$ and cannot be recovered to $C$ once the traffic system becomes congested.

When $(1-\Delta)C< d^-\leq C$, we also set $u(t)=v_1$ and  rewrite \refe{system_8au} as
\bqn
\der{}{t}k&=&\frac 1{l_0}(\min\{d^-, w(k_j-k)\}-\min\{v_f k, (1-\Delta \cdot I_{k>k_1})C\}). \label{syssta2}
\eqn
\ben
\item For the uncongested equilibrium state at $k^*=\frac{d^-}{v_f}$ and $g^*=d^-$, we denote a small disturbance in traffic density by $z=k-k^*$, then \refe{syssta2} is equivalent to
\bqs
\der{}{t}z&=&-\frac {v_f}{l_0}z. 
\eqs
Therefore $z=0$ is an exponentially stable equilibrium point, and $k^*$ is exponentially stable. However, if a disturbance is sufficiently large such that $k>k_1$, then \refe{syssta2} becomes
\bqs
\der{}{t}k&=&\frac 1{l_0}(d^--(1-\Delta)C)>0,
\eqs
and the system will further drift away from $k^*$.

\item For the congested equilibrium state at $k^*=k_2$ and $g^*=(1-\Delta)C$, we denote a small disturbance in traffic density by $z=k-k_2$, then \refe{syssta2} is equivalent to
\bqs
\der{}{t}z&=&-\frac w{l_0} z.
\eqs
Therefore $z=0$ is an exponentially stable equilibrium point, and $k^*=k_2$ is exponentially stable.
\een
Thus when $(1-\Delta)C< d^-\leq C$, with $u^*=v_1$, any large positive disturbance in density would drive the system away from the optimal equilibrium state, and the discharging flow-rate will be reduced to $(1-\Delta)C$ and cannot be recovered to $C$ once the traffic system becomes congested.

In summary, we can see that $u^*=v_1$ is an optimal speed limit when $(1-\Delta)C< d^-\leq C$. However, when there exists a sufficiently large disturbance, the open-loop control system will drift away from the optimal uncongested equilibrium state and converges to the exponentially stable, congested equilibrium state, in which the discharging flow-rate is the dropped capacity. Therefore, the constant optimal speed limit $u^*=v_1$ is not an effective control strategy.

\section{A closed-loop control system with a PI feedback controller}
In a state feedback control system, the control signal is also a function of $k$, i.e., 
\begin{equation}
u = u(k).
\end{equation} 
As the result, \refe{controlsys} becomes an autonomous system, and its controller design can be simplified.
In addition, the speed limit has to be within the following range
\bqn
u_{min}\leq u(t) \leq v_f,
\eqn
where $u_{min}<v_2$ is the minimum speed limit. A trivial choice is $u_{min}=0$.
In reality, however, a smaller $u$ leads to a longer distance for vehicles to accelerate away from the VSL zone. This would require a longer zone 0 and can cause interference with further upstream off-ramp flows. Therefore, $u$ should be as large as possible. Another constraint on $u$ is that it should set as discrete values at, for examples, 30 mph, 35 mph, 40 mph, 45 mph, and so on, so that drivers can easily follow the speed limits. But in this study we ignore the discrete nature of the speed limit and assume that $u(t)$ is continuous. 

From the analysis in Section 3.2, we find that $v_1$ is always an optimal speed limit. Therefore, it is reasonable to set the target density of $k(t)$ to be $\bar k=k_1$, since both of them correspond to a discharging flow-rate of $C$. 
In this section, we consider the following PID (proportional-integral-derivative) controller, which is widely used in industry \citep[][Chapter 10]{astrom2008feedback}:
\bsq\label{pidcontroller}
\bqn
u(t)&=&v_1+ \alpha e(t)+ \beta\int_{\tau=0}^t e(\tau) d\tau+\gamma \der{e(t)}t,
\eqn 
where $e(t)=k_1-k(t)$ is the error signal. 
\esq

Usually $\alpha$, $\beta$, and $\gamma$ are non-negative parameters to be determined. 
In the controller, \refe{pidcontroller}, the proportional term $\alpha e (t)$, which is also called the ``P-term'', generates control signals proportional to the error of $k$. The integral term $\beta \int_0^t e(\tau) d \tau$, also caller the ``I-term'', guarantees that $k(t)$ converge to the equilibrium state, even subject to uncertainties in the parameters of \refe{controlsys}. 
Usually the derivative term is  used to improve the transient performance of a control system, but in this study we focus on the equilibrium states and their stability. Thus we set $\gamma=0$, and \refe{pidcontroller} is effectively a PI controller. A detailed introduction of PI and PID controllers can be found in \citep[][Chapter 10]{astrom2008feedback}.

Thus the traffic system, \refe{controlsys}, and the PI controller, \refe{pidcontroller}, form a closed-loop feedback control system, for which we will analyze the equilibrium states and their stability in the remainder of this section.

\subsection{Equilibrium states}
We denote the asymptotic equilibrium states of the closed system, \refe{controlsys} and \refe{pidcontroller}, by $k(t)=k^*$, $u(t)=u^*$, and $e(t)=e^*$. Thus at a very large $t$ we have
\bsq \label{pidequ}
\bqn
&&\min\{d^-, \frac{u^*}{u^*+w} wk_j, w(k_j-k^*)\}=\min\{v_f k^*, C(1-\Delta \cdot I_{k^*>k_1})\},\\
&&e^*=k_1-k^*,\\
&&u^*=v_1+\alpha e^*+\beta e^* t, \label{pideqn3}\\
&&u_{min}\leq u^*\leq v_f.
\eqn
\esq

We require that at least one of the P- and I-terms be used; i.e., $\alpha>0$ or $\beta>0$. Otherwise, $\alpha=\beta=0$, \refe{pidequ} is equivalent to \refe{equpoint}, and the system can have multiple equilibrium states, including a congested one.  

When $d^-\geq C$, we can have the following equilibrium states
\ben
\item It is easy to check that \refe{pidequ} always has the following uncongested equilibrium state: $k^*=k_1$, $u^*=v_1$, and $e^*=0$. 
\item When $\beta=0$, if $\alpha\geq \frac{v_1-v_2}{k_2-k_1}$, \refe{pidequ} can have another equilibrium state: $u^*=v_2$, $e^*=-\frac{v_1-v_2}{\alpha}<0$, and $k^*=k_1+\frac{v_1-v_2}{\alpha}\leq k_2$. The corresponding discharging flow-rate is $g^*=(1-\Delta) C$.
\item When $\beta=0$, if $\alpha\leq \frac{v_1-v_2}{k_2-k_1}$, \refe{pidequ} can have another equilibrium state: $k^*=k_2$, $e^*=k_1-k_2<0$, and $u^*=v_1 -\alpha (k_2-k_1)\geq v_2$. The corresponding discharging flow-rate is $g^*=(1-\Delta) C$.
\item  It is no longer an equilibrium state with $k^*=k_2$ and $u^*=v_1$, since $e^*<0$ and \refe{pideqn3} is violated. 
\item When $\beta>0$, there is only one equilibrium state.
\een
We can see that the I-term successfully removes the congested equilibrium state with $g^*=(1-\Delta)C$ in the open system. Therefore, the I-term has to be included; i.e., we require $\beta>0$ in \refe{pidcontroller}.  

When $d^-< C$, we can see that \refe{pidequ} always has only one uncongested equilibrium state: $k^*=\frac{d^-}{v_f}$, $u^*=v_f$, and $e^*=k_1-k^*>0$. Since $\beta>0$, there is no other equilibrium states.
Therefore, with $\beta>0$, the PI controller in \refe{pidcontroller} will lead to ideal equilibrium states for all demand patterns.

\subsection{Stability of equilibrium states}
In this subsection, we analyze the stability of the only uncongested equilibrium state at $k^*=\min\{\frac{d^-}{v_f}, k_1\}$, $u^*=\cas{{ll} v_1, & d^-\geq C\\ v_f, & d^-<C}$, and $e^*=k_1-k^*$. 
In this case, 
the closed-loop system, \refe{controlsys} and \refe{pidcontroller}, is equivalent to
\bsq \label{pisys}
\bqn
\der{}{t}k(t)&=&\frac 1{l_0}( \min\{d^-, \frac{u(t)}{u(t)+w} wk_j, w(k_j-k(t))\}-\min\{v_f k(t), C(1-\Delta \cdot I_{k(t)>k_1})\}),\\
\der{}{t}u(t)&=&-\alpha \der{k(t)}{t}+\beta (k_1-k(t)).
\eqn
\esq

When $d^-\geq C$, we denote the disturbances by $\epsilon(t)=u(t)-v_1$ and $z(t)=k(t)-k_1$.In the neighborhood of $(k_1, v_1)$, both $z(t)$ and $\epsilon(t)$ are small, then  $\frac{u(t)}{u(t)+w} wk_j\leq \min\{d^-, w(k_j-k(t))\}$, and system \refe{pisys} can be linearized as
\bsq \label{smallosc}
\bqn
\der{}{t} z(t)&\approx&\frac 1{l_0}(  k_3 \epsilon (t)-\min\{v_f z(t), -C \Delta \cdot I_{z(t)>0}\}),\\
\der{}{t} \epsilon (t)&=&-\alpha \der{z(t)}{t}-\beta z(t),
\eqn
\esq
where $k_3=\frac{w^2 k_j}{(v_1+w)^2}$.
Thus when $z(t)\leq 0$, \refe{smallosc} is equivalent to
\bsq\label{smallosc1n}
\bqn
\der{}{t} z(t)&=&-\frac{v_f}{l_0} z(t)+\frac {k_3}{l_0} \epsilon (t),\\
\der{}{t} \epsilon (t)&=& (\frac{\alpha v_f}{l_0}-\beta) z(t)-\frac{\alpha k_3}{l_0} \epsilon(t).
\eqn
\esq
When $z(t)>0$, \refe{smallosc} is equivalent to
\bsq\label{smallosc1p}
\bqn
\der{}{t} z(t)&=&\frac {k_3}{l_0} \epsilon (t)+\frac{C \Delta}{l_0},\\
\der{}{t} \epsilon (t)&=& -\beta z(t)-\frac{\alpha k_3}{l_0} \epsilon(t)-\frac{\alpha C \Delta}{l_0}.
\eqn
\esq

The linearized system, \refe{smallosc1n} and \refe{smallosc1p}, is a switched linear dynamical system, for which the equilibrium point of \refe{smallosc1n} is the origin. But the equilibrium point of \refe{smallosc1p} is at $z^*=0$ and $\epsilon^*=-\frac{C\Delta}{k_3}$. Thus \refe{smallosc} and \refe{pisys} can reach the origin only from the left half of the $z-\epsilon$ space when $z\leq 0$. 
Therefore, the stability of the closed-loop traffic system is highly related to that of  the switched dynamical system showed that \refe{smallosc1n} and \refe{smallosc1p}. However, a theoretical analysis of this system is beyond the scope of this study \citep{sun2011stability}. In the following section, we will present some simulation-based studies. 

When $d^-<C$, we denote the disturbances by $\epsilon(t)=u(t)-v_f$ and $z(t)=k(t)-\frac{d^-}{v_f}$. Assuming that both $\epsilon(t)$ and $z(t)$ are small, then  $d^-<\min\{\frac{u(t)}{u(t)+w} wk_j, w(k_j-k(t))\}$, and system \refe{pisys} can be linearized as
\bsq
\bqn
\der{}{t} z(t)&\approx&-\frac {v_f}{l_0} z,\\
\der{}{t} \epsilon (t)&=&(\frac {\alpha v_f}{l_0}-\beta) z(t).
\eqn
\esq
Therefore the equilibrium state is always asymptotically stable whether $\alpha=0$ or not.

\section{Numerical simulations of the closed-loop control system}
In this section, we consider the VSL control problem for a lane-drop bottleneck from two lanes to one lane. The length of zone 0 is in the order of 700 meters \citep{carlson2010optimal}. Here we set the length of zone 0 to be ${l_0}=600$ m. We pick the following parameters \citep{yang2011calibration}: $v_f=$30 m/s, $w=35/8$ m/s, and $k_j=2/7$ veh/m. Then $k_c=2/55$ veh/m. 
In this case, lane-changes do not reduce the downstream link's capacity, which is $C=\frac 12 v_f k_c=6/11$ veh/s \citep{jin2013multi}. We assume that the capacity drop magnitude is $\Delta=20\%$.  Then we can calculate $k_0$, $k_1$, $k_2$, $k_3$, $v_1$, and $v_2$ correspondingly. We set the target density $\bar k=k_1$ and $u_{min}=0.5$ m/s. Here we numerically solve \refe{pisys} by using Euler's method with a time-step size of 1 s

\subsection{Stability property with I-controller and PI-controller}
Let $d^-=2 C$ and $k(0)=2 k_1$. From the analyses in the preceding sections, we know that the system reaches the congested equilibrium state at $k^*=k_2$ without VSL control or with an open loop control $u^*=v_1$. 

When we apply an I-controller with $\beta=4$, the simulated density, speed limit, and discharging flow-rate are shown in \reff{I-result}(a). From the figure, we can see that the system converges to an asymptotically periodic state, with an average discharging flow-rate of $C$, which is the maximum discharging flow-rate. In this case, the switched dynamical system, \refe{smallosc1n} and \refe{smallosc1p}, exponentially converges to the origin, and the capacity drop is completely avoided. However, when we apply an I-controller with $\beta=20$, the simulated density, speed limit, and discharging flow-rate are shown in \reff{I-result}(b). From the figure, we can see that the system converges to an asymptotically periodic state, with an average discharging flow-rate of $0.7988 C$, which smaller than $(1-\Delta)C$. This confirms the existence of a limit cycle in the switched dynamical system, \refe{smallosc1n} and \refe{smallosc1p}.

\bfg \bc
$\ba{c@{\hspace{0in}}c}
\includegraphics[height=2.5in]{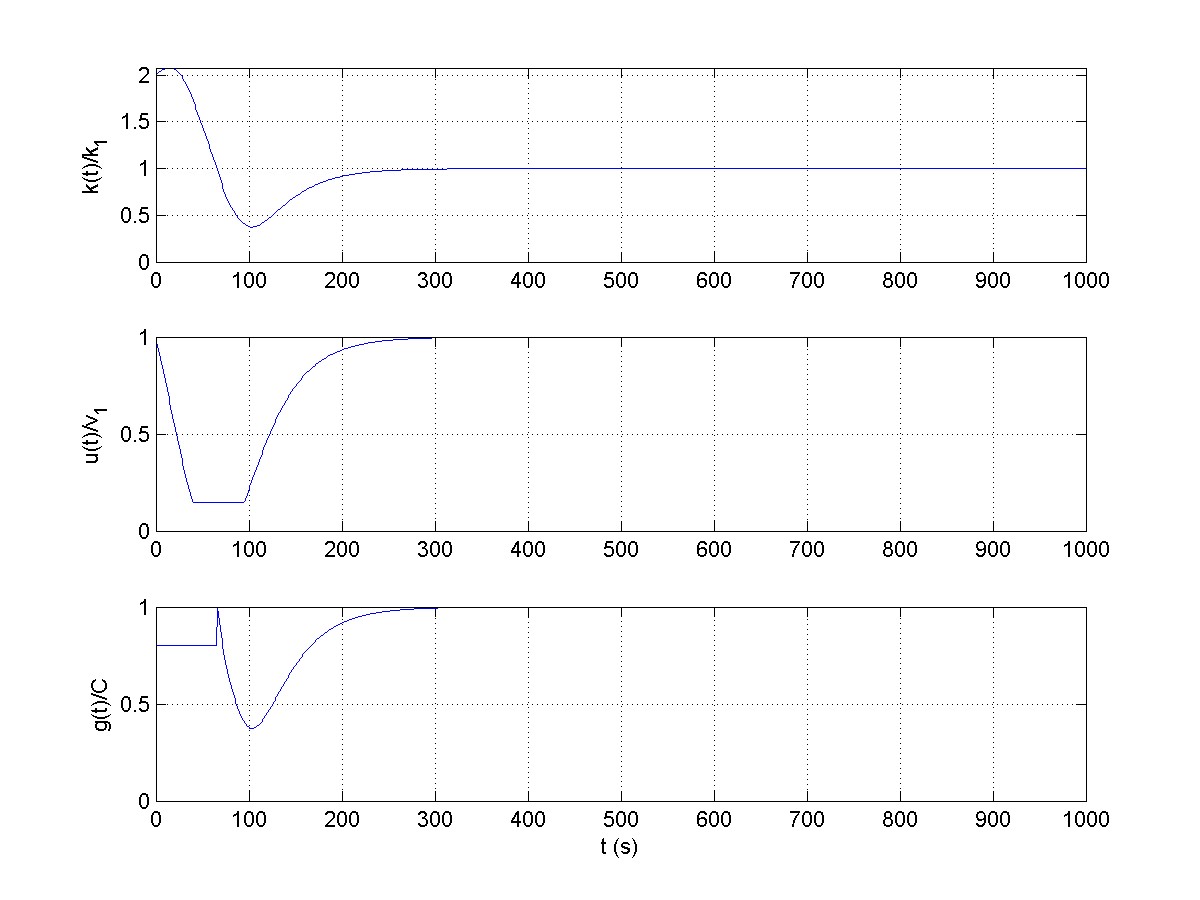} &
\includegraphics[height=2.5in]{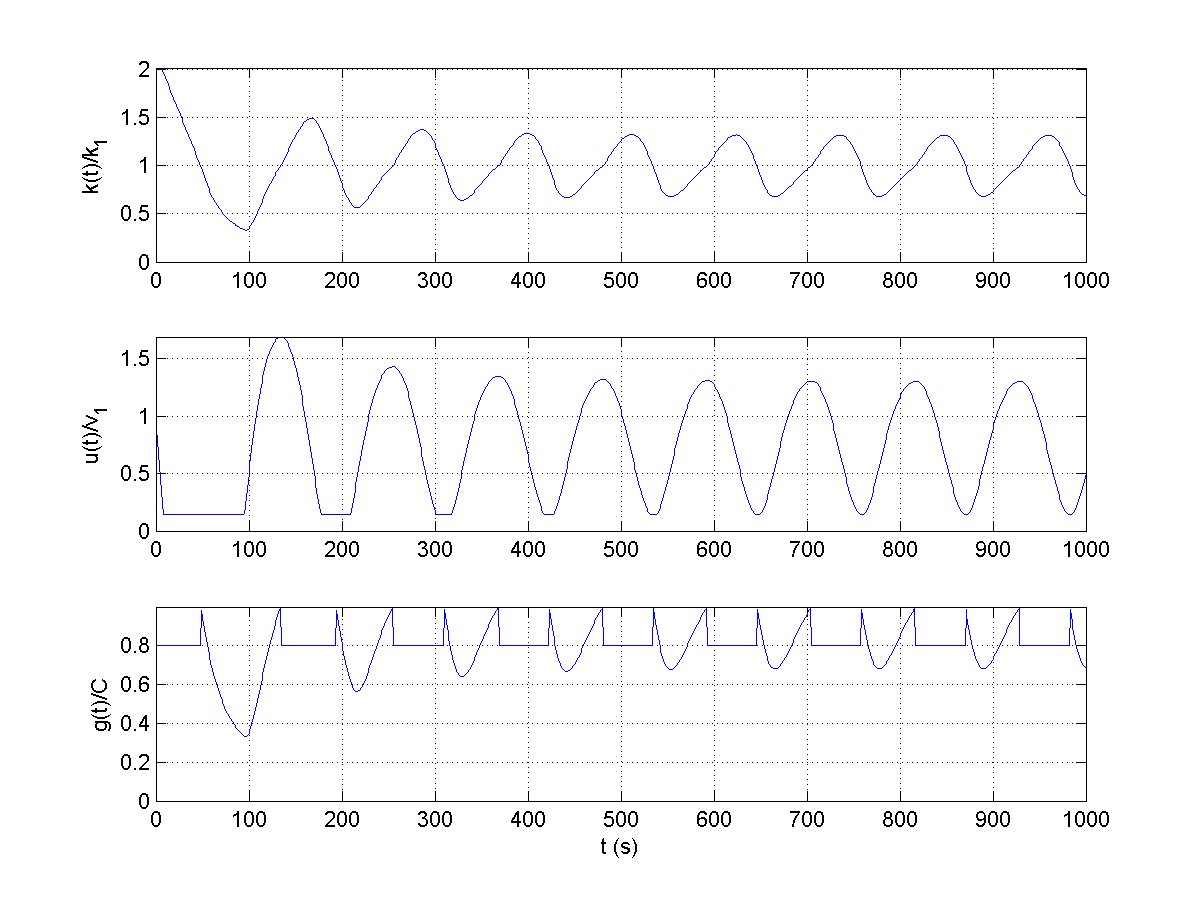} \\
\multicolumn{1}{c}{\mbox{\bf (a)}} &
    \multicolumn{1}{c}{\mbox{\bf (b)}}
\ea$
\caption{Simulation results of \refe{pisys} with an I-controller: (a) $\beta=4$; (b) $\beta=20$} \label{I-result} \ec 
\efg

When we apply a PI-controller with $\alpha=400$ and $\beta=20$, the simulated density, speed limit, and discharging flow-rate are shown in \reff{PI-result}(a). From the figure, we can see that the system converges to an asymptotically periodic state, with an average discharging flow-rate of $0.9202 C$, which is greater than $(1-\Delta)C$ but smaller than C. Such a PI-controller clearly outperforms the I-controller, but the switched dynamical system, \refe{smallosc1n} and \refe{smallosc1p}, still converges to a limit cycle.
In contrast, when we apply a PI-controller with $\alpha=500$ and $\beta=20$, the simulated density, speed limit, and discharging flow-rate are shown in \reff{PI-result}(b). From the figure, we can see that the system converges to an asymptotically periodic state, with an average discharging flow-rate of $C$, which is the maximum discharging flow-rate. In this case, the switched dynamical system, \refe{smallosc1n} and \refe{smallosc1p}, exponentially converges to the origin, and the capacity drop is completely avoided.

\bfg \bc
$\ba{c@{\hspace{0in}}c}
\includegraphics[height=2.5in]{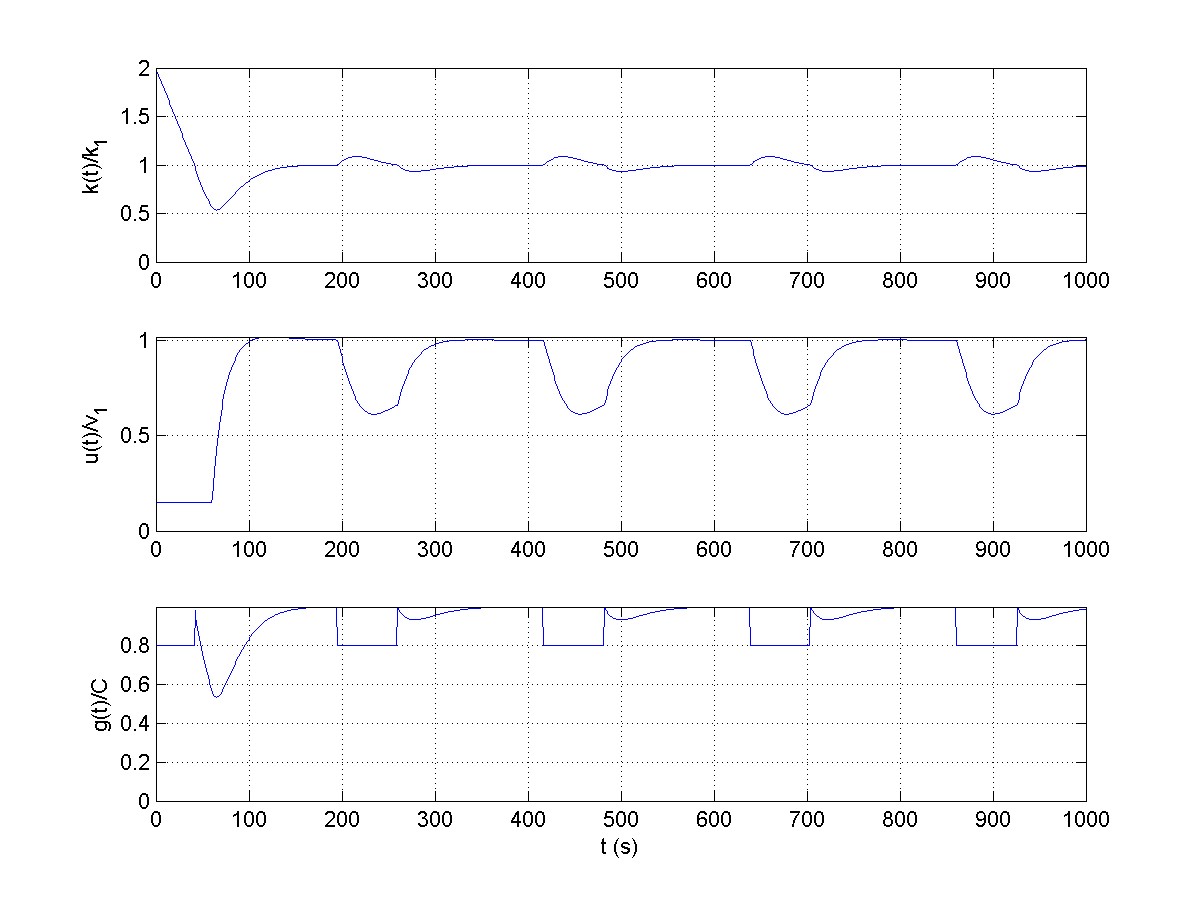} &
\includegraphics[height=2.5in]{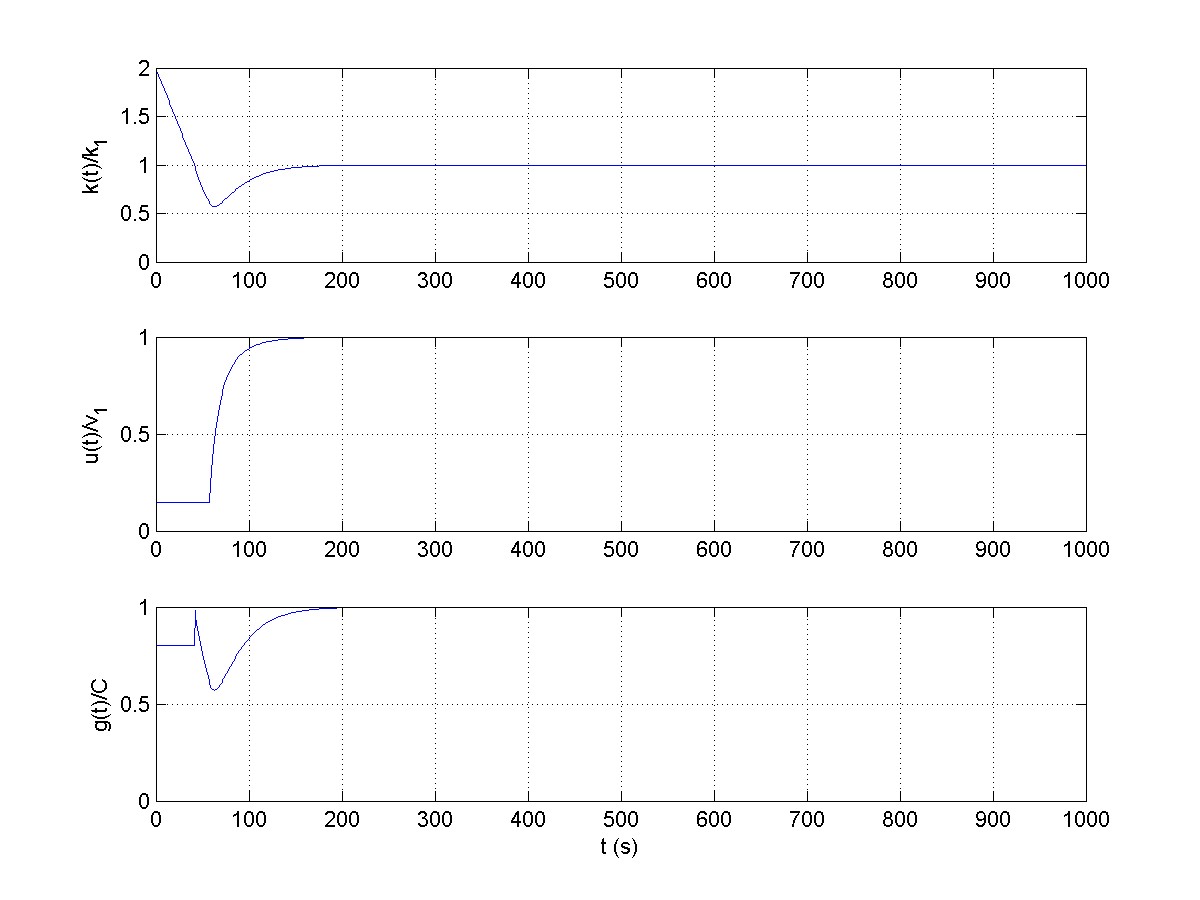} \\
\multicolumn{1}{c}{\mbox{\bf (a)}} &
    \multicolumn{1}{c}{\mbox{\bf (b)}}
\ea$
\caption{Simulation results of \refe{pisys} with a PI-controller: (a) $\alpha=400$ and $\beta=20$; (b) $\alpha=500$ and $\beta=20$} \label{PI-result} \ec 
\efg

The simulation results confirm that the closed-loop system with both I- and PI-controllers can be stable and effective, if the coefficients are appropriately chosen for different demand patterns, $d^-$, initial conditions, $k(0)$, the length of zone 0, $l_0$, and the minimum speed limit, $u_{min}$. In particular, in the following studies, we only consider the I-controller with $\beta=4$ and the PI-controller with $\alpha=500$ and $\beta=20$. 

\subsection{Robustness with respect to the estimation of $k_1$}
In reality, since both $C$ and $v_f$ can vary with traffic conditions, the target density, $k_1=\frac{C}{v_f}$, cannot be accurately measured. That is, the target density can be under- or over-estimated. In this subsection, we introduce an error in the target density:
\bqn
\bar k&=&(1+\xi) k_1. \label{xidef}
\eqn
Let $d^-=2 C$,  $k(0)=2 k_1$. In this subsection we study the impacts of $\xi$ on the performance of the control system with an I-controller ($\beta=4$). Note that an error of $\xi k_1$ in estimating $k_1$ equals an error of $-\xi k_1$ in measuring the state variable $k$.

\begin{figure} \bc
$\ba{c@{\hspace{0in}}c}
\includegraphics[height=2.5in]{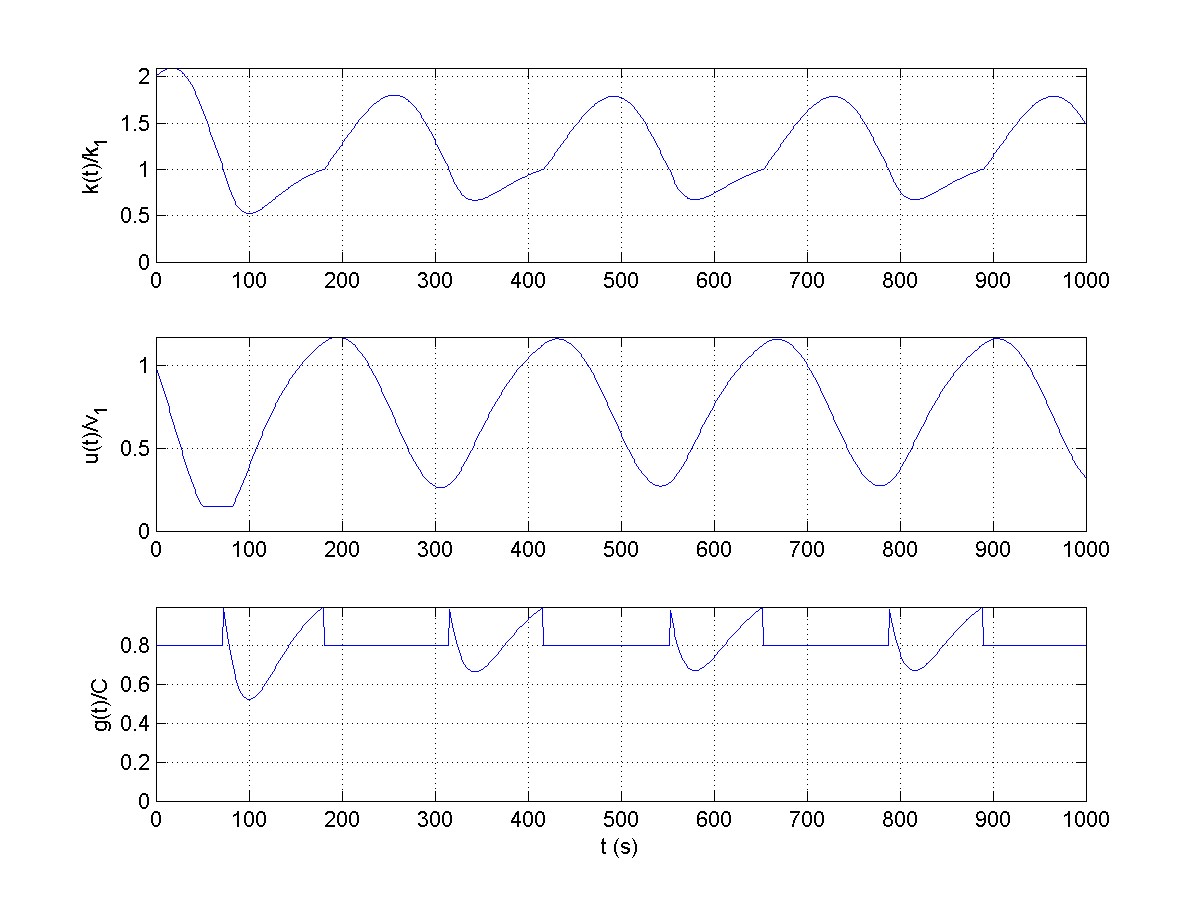} &
\includegraphics[height=2.5in]{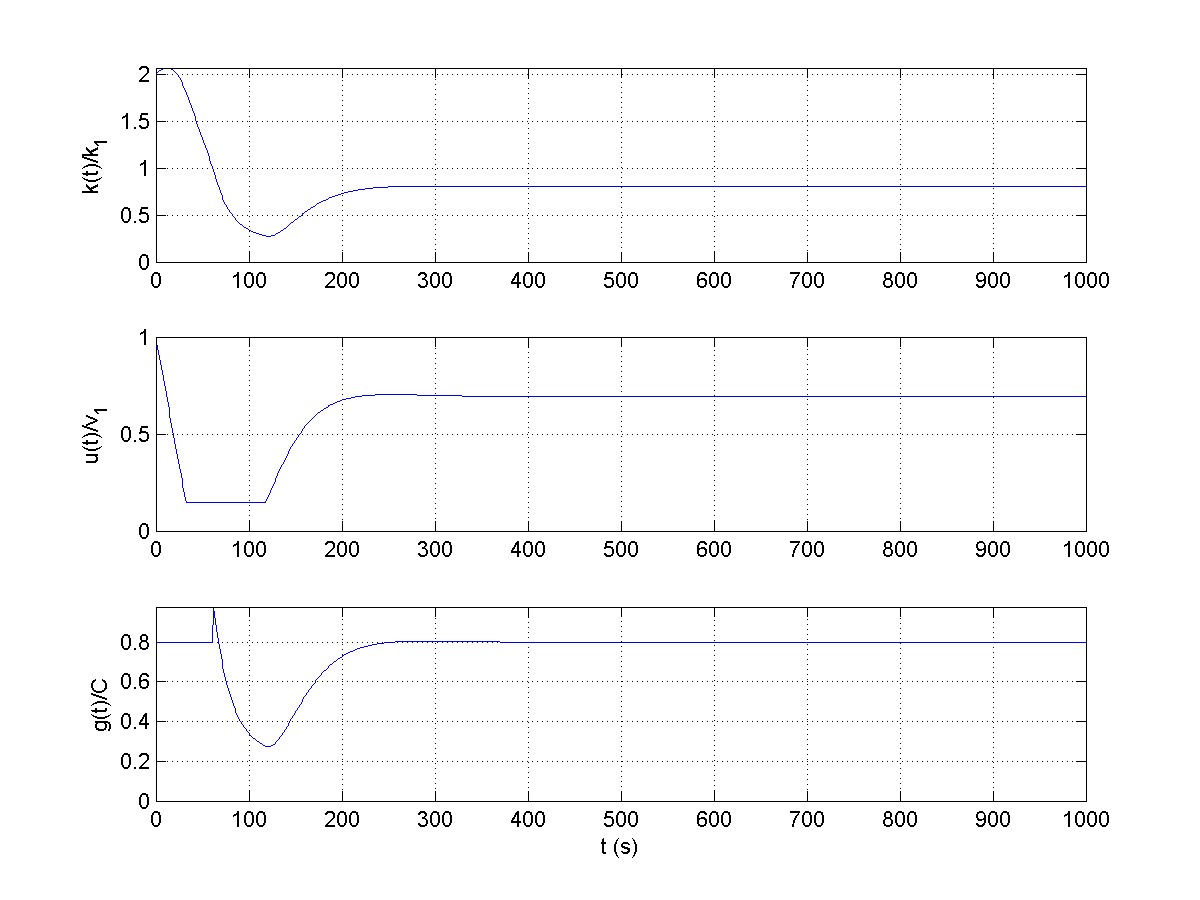} \\
\multicolumn{1}{c}{\mbox{\bf (a)}} &
    \multicolumn{1}{c}{\mbox{\bf (b)}}
\ea$
\caption{Simulation results of \refe{pisys} with inaccurate target density: (a) $\bar k=1.1 k_1$; (b) $\bar k=0.9 k_1$} \label{targetdensity} \ec 
\end{figure}

In \reff{targetdensity}, we demonstrate the simulation results of density, speed limit, and discharging flow-rate for an over-estimated target density $\bar k=1.1 k_1$ and an under-estimated target density $\bar k=0.9 k_1$. From the simulation results, we can see that the traffic system converges to a periodic state when the target density is over-estimated, as shown in \reff{targetdensity}(a), and to an uncongested equilibrium state when the target density is under-estimated. In particular, when the target density is under-estimated, the equilibrium state is at $k^*=\bar k<k_1$ and $g^*=v_f \bar k$, and capacity drop does not occur; but when the target density is over-estimated, the I-controller still helps to improve the performance of the system since the average discharging flow-rate $g^*=0.81 C>(1-\Delta)C$, but capacity drop occurs periodically.

In \reff{gstar_xi}, we further demonstrate the asymptotic average discharging flow-rate $g^*$ for different $\xi$. From the figure, we can clearly see that $g^*$ is a discontinuous function in $\xi$,  periodic capacity drop occurs when the target density is over-estimated, but the discharging flow-rate continuously changes with $\xi$ when the target density is under-estimated. Therefore, when implementing the VSL control, we would rather use a smaller target density than a bigger one.
However, how to adaptively change the target density is beyond the scope of this research.

\bfg\bc
\includegraphics[width=4in]{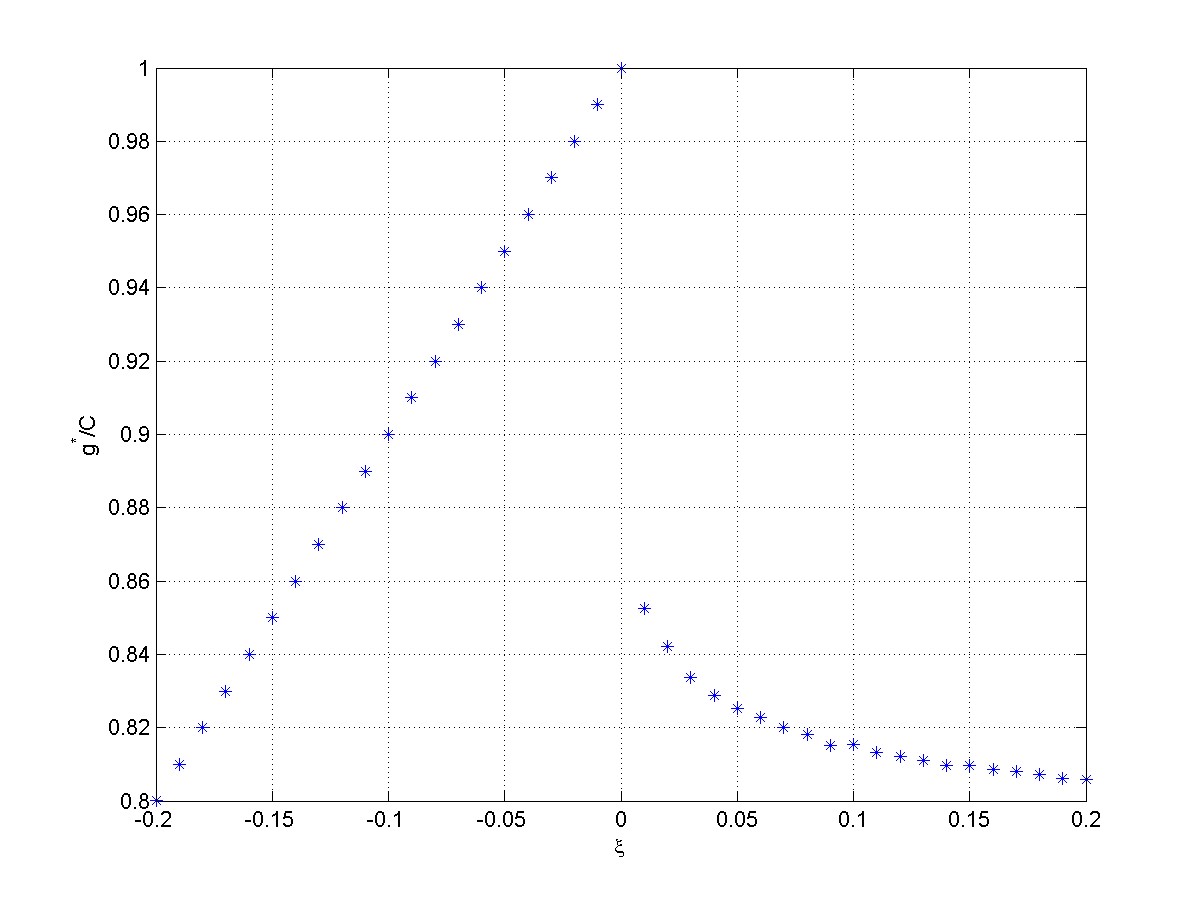}
\caption{Equilibrium discharging flow-rate versus $\xi$ in \refe{xidef}}\label{gstar_xi}
\ec\efg

\subsection{Robustness with respect to variations in the demand level} 
In this subsection, we consider dynamic arrival patterns for the whole road in \reff{lanedrop_vsl_ill} with random noises as follows ($t\in[0,8000]$)
\bqn
r(t)&=&\max\{0, C \min\{1, 0.0005 t, 1-0.0005(t-4000)\}+\tilde r\},
\eqn
where the random noise $\tilde r$ follows a normal distribution $N(0,0.02 C)$.
We model the VSL zone and the upstream queueing zone as a point queue:
\bsq\label{pointqueue}
\bqn
\der{\lambda(t)}{t}&=&r(t)-f(t),\\
d^-(t)&=&\cas{{ll} v_f k_c, & \lambda(t)>0\\ \min\{v_f k_c, r(t)\}, & \lambda(t)=0},
\eqn
\esq
where $\lambda(t)\geq 0$ is the queue size.
In the discrete form, \refe{pointqueue} can be written as
\bsq\label{discretepq}
\bqn
\lambda(t+\dt)&=&\lambda(t)+\dt(r(t)-f(t)),\\
d^-(t)&=&\min\{v_f k_c,\frac{\lambda(t)}{\dt}+r(t)\}.
\eqn
\esq
Since $f(t)\leq d^-(t)$, it is guaranteed that $\la(t+\dt)\geq 0$ if $\la(t)\geq 0$.
We assume that the target density $\bar k=k_1$ is accurately estimated. We let $\la(0)=0$, and $k(0)=0$. 

First we consider the I-controller with $\beta=4$. 
The simulation results are given in \reff{vsl_pid_random}, where all solid lines are for the results with the PI controller in \refe{pidcontroller}, and all dashed lines for those without VSL control; i.e., $u(t)=v_f$. From the figures, we can see that, when the arrival flow-rate is less than $C$ before $t=2000$ s, the queue size is 0, the in-flux equals the arrival flow-rate, as predicted by \refe{pointqueue}, and the traffic system has the same performance with or without control. However, at $t=2000$ s, the arrival flow-rate is greater than $C$ due to the disturbance, the PI controller starts to reduce the speed limit in the VSL zone towards the optimal speed limit $v_1$, as shown in \reff{vsl_pid_random}(f), and the traffic density in zone 0 is controlled toward $k_1$, as shown in \reff{vsl_pid_random}(e). After that, the VSL control successfully increases the in-flux, as shown in \reff{vsl_pid_random}(c), prevents the occurrence of capacity drop and maintains a high discharging flow-rate, as shown in \reff{vsl_pid_random}(d), and maintains a small queue in the VSL and upstream region, as shown in \reff{vsl_pid_random}(b). Thus we can see that the PI controller is robust with respect to random variations in the demand patterns.

\bfg\bc
\includegraphics[width=5in]{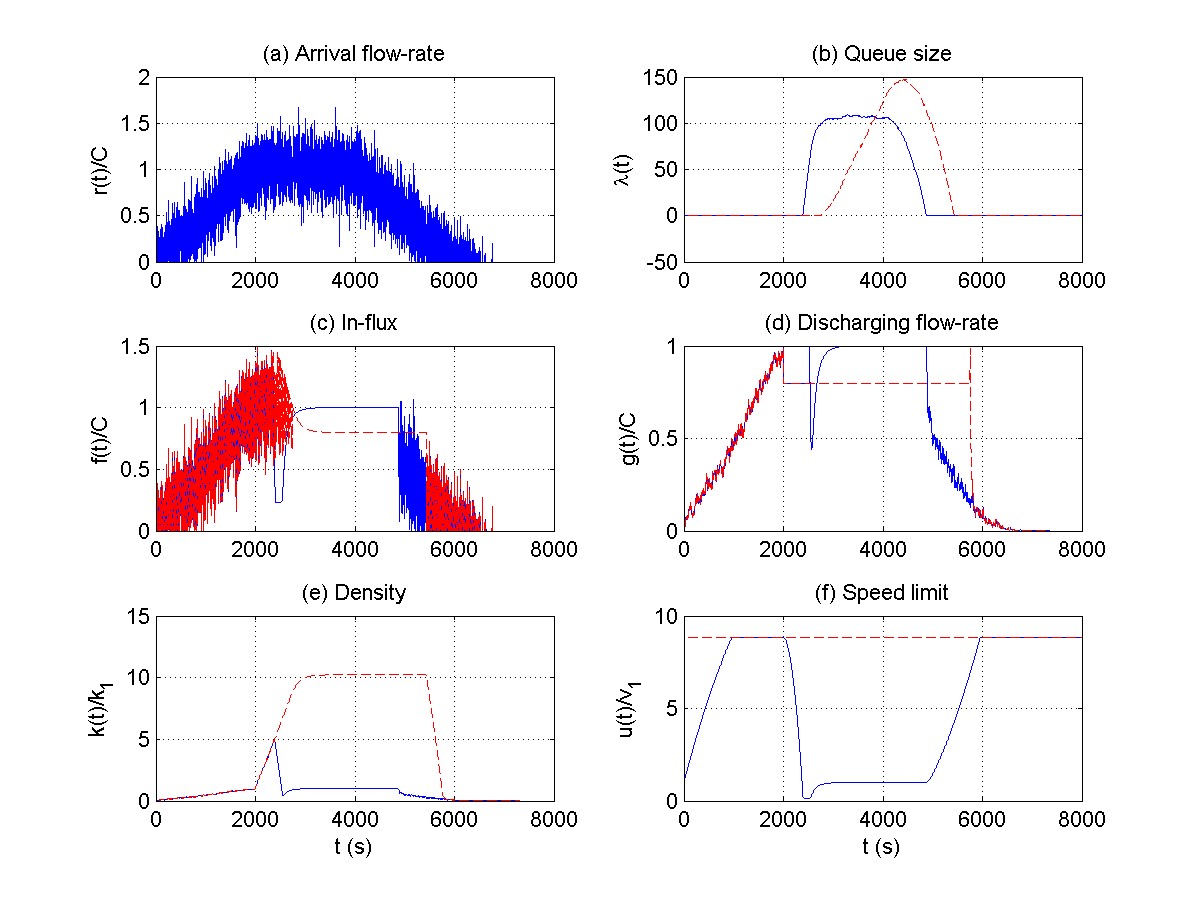}
\caption{Simulation of traffic dynamics in a lane-drop region with (solid lines) and without (dashed lines) the VSL control}\label{vsl_pid_random}
\ec\efg

In \reff{vsl_pid_random_cumulative}, we demonstrate the cumulative arrival flows, $\int_0^t r(\tau) d\tau$, and departure flows, $\int_0^t g(\tau) d\tau$, with or without VSL control.
In total, 2223 vehicles travel through the lane-drop region during 8000 s. The average travel time through the bottleneck is 122 seconds with the VSL control and 268 seconds without VSL control. Thus in this case the VSL control is able to save 55\% of the travel time. 

\bfg\bc
\includegraphics[width=5in]{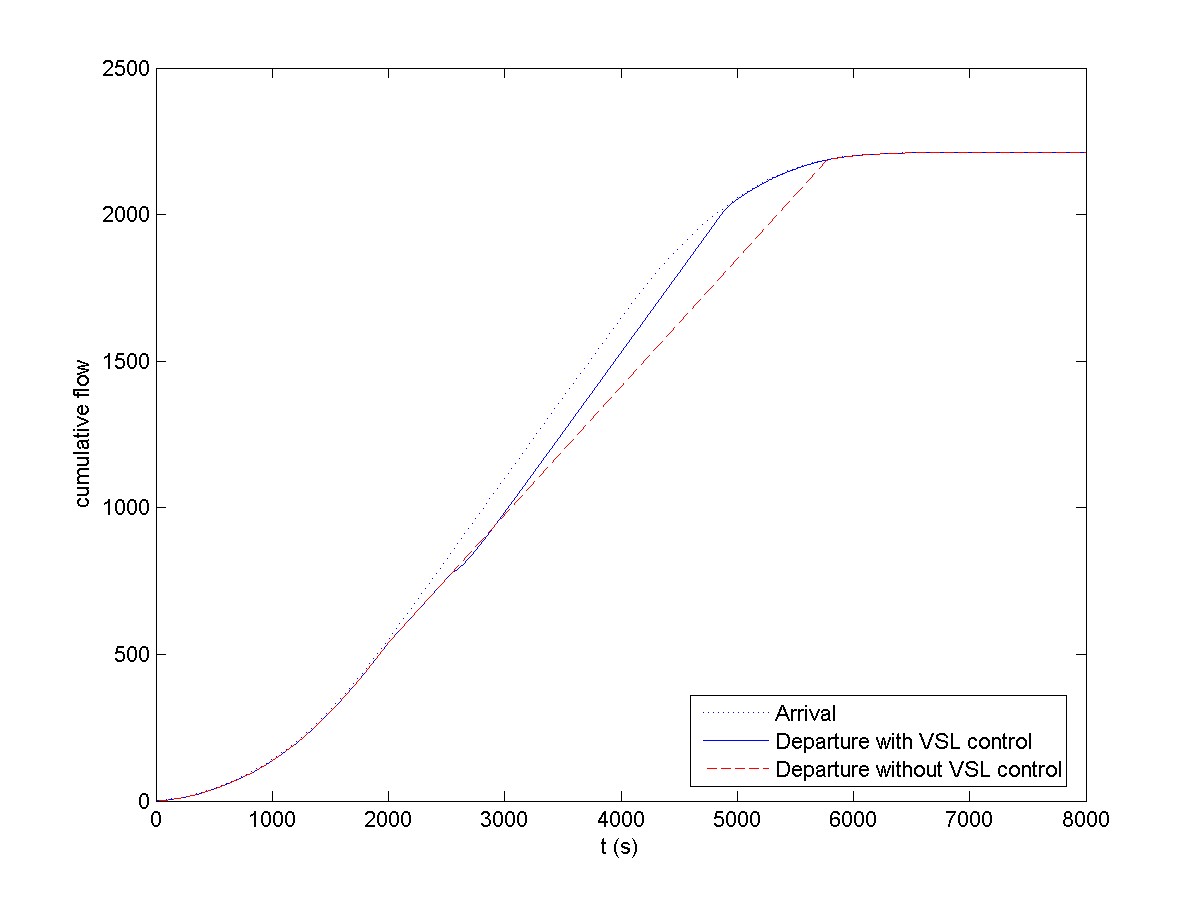}
\caption{Cumulative flows with (solid lines) and without (dashed lines) VSL control}\label{vsl_pid_random_cumulative}
\ec\efg

Further in \reff{vsl_pid_random_save}, we show the reduction ratio in travel time, defined as \[1-\frac{\m{Travel time with VSL control}}{\m{Travel time without control}}.\]
First, we confirm that the VSL control cannot improve the system's performance without capacity drop ($\Delta=0$). Second, the reduction ratio increases with the capacity drop magnitude in a nonlinear fashion. Third, the exact reduction ratio depends on the road geometry, traffic conditions, and the parameters in the PI controller.

\bfg\bc
\includegraphics[width=5in]{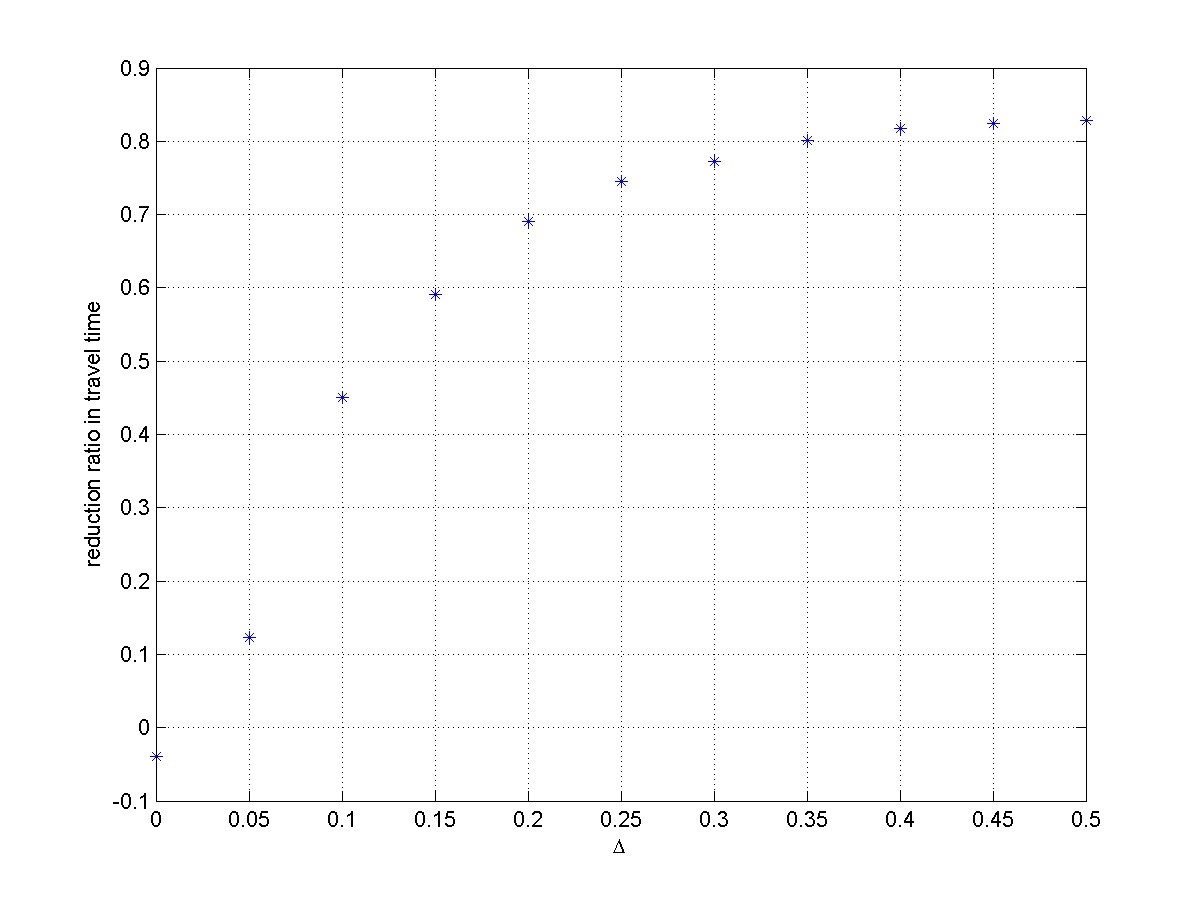}
\caption{The reduction ratio in travel time due to the VSL control}\label{vsl_pid_random_save}
\ec\efg

In \reff{vsl_pid_random_pi}, we show the simulation results of a PI controller with $\alpha=500$ and $\beta=20$. Compared with \reff{vsl_pid_random}, we can see that the PI controller is also effective and robust with respect to random variations in the demand patterns. 
\bfg\bc
\includegraphics[width=5in]{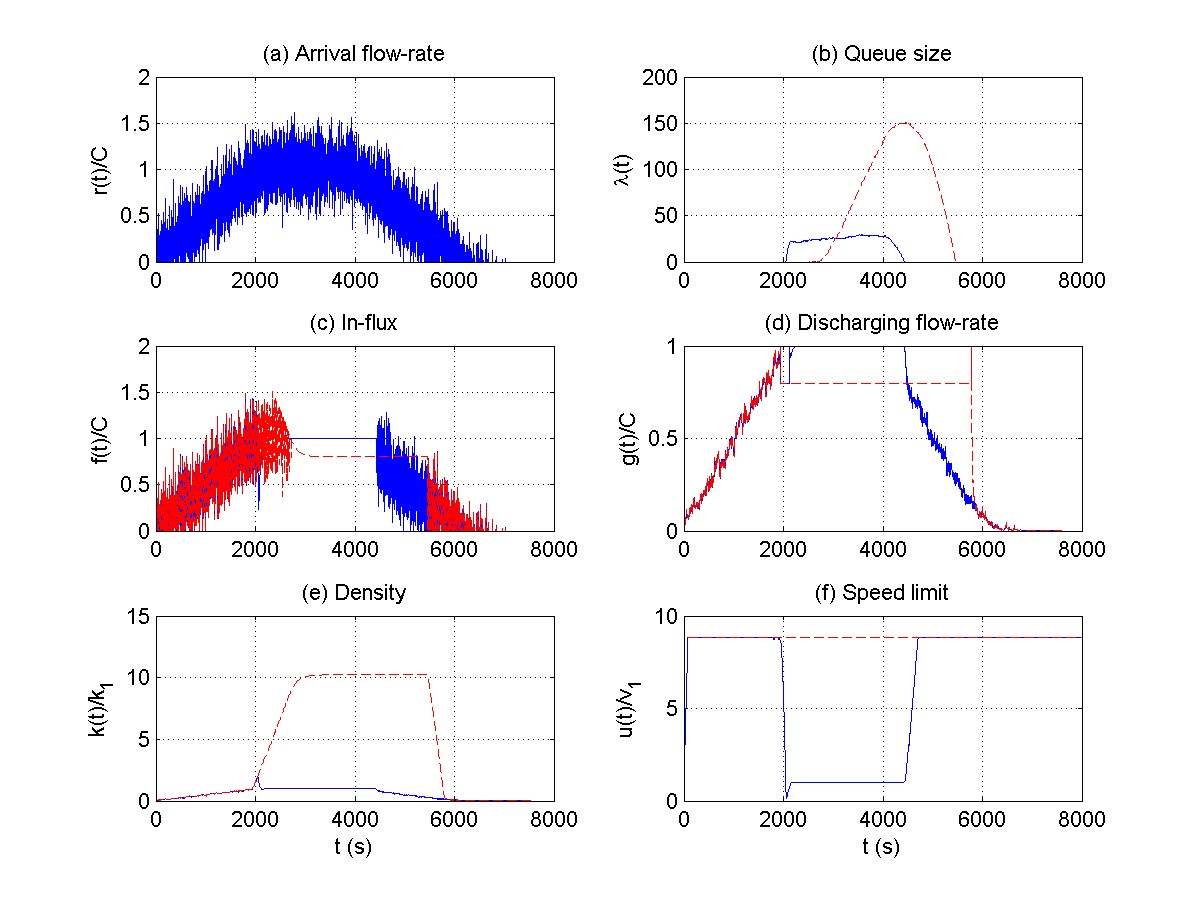}
\caption{Simulation of traffic dynamics in a lane-drop region with (solid lines) and without (dashed lines) the VSL control: $\alpha=500$ and $\beta=20$}\label{vsl_pid_random_pi}
\ec\efg

\section{A feedback control system based on the LWR model}
In this section, we apply the PI controller developed in the preceding sections to control the LWR model \refe{lwr}. We first discretize the model with the entropy condition, \refe{entropycond}, into  the Cell Transmission Model \citep{daganzo1995ctm}, in which there are $n$ cells with the cell length of $\dx$ and a corresponding time-step size $\dt$:
\bsq \label{ctm}
\bqn
\r_i^{j+1}&=&\r_i^j+\frac{\dt}{\dx} (q_{i-\frac 12}^j -q_{i+\frac12}^j), \quad i=1,\cdots,n\\
\delta_{i}^j&=&D(\r_i^j), \quad \sigma_i^j=S(\r_i^j), \quad i=1,\cdots,n\\
q_{i-\frac 12}^j&=&\min\{\delta_{i-1}^j,\sigma_i^j\}, \quad i=2,\cdots,n \label{discreteflux}
\eqn
where the CFL number $v_f \frac\dt\dx \leq 1$ \citep{courant1928CFL}, $\r_i^{j}$ is the average density in cell $i$ at $j\dt$, $q_{i-\frac12}^j$ is the flux at $x=(i-1) \dx$ during $[j\dt,(j+1)\dt]$, and \refe{discreteflux} is the discrete version of \refe{entropycond}.
At the upstream and downstream boundaries, we apply \refe{infloweqn-lwr} and \refe{dischargingflow-lwr} to calculate the in-flux, $q_{\frac 12}^j=f^j$, and the out-flux, $q_{n+\frac12}=g^j$ as follows:
\bqn
q_{\frac 12}^j&=&\min\{d^j, \frac{u^j}{u^j+w} wk_j, w(k_j-\r_1^j)\},\\
q_{n+\frac12}&=&\min\{v_f \r_n^j, C(1-\Delta \cdot I_{\r_n^j>k_1})\}, \label{ctmdischarging}
\eqn
\esq
where $d^j=d^-(j\dt)$ is the traffic demand, and $u^j$ the speed limit at $j\dt$.

Since capacity drop is directly triggered by $\r_n^j$ in \refe{ctmdischarging}, we feedback the traffic density in the last cell, $\r_n^j$, to determine the speed limit in the PI controller, \refe{pidcontroller}. We still set the target density to be $k_1$. Then the speed limit can be updated by
\bqn
u^{j+1}&=&u^j-\alpha (\r_n^{j+1}-\r_n^j)+\beta  (k_1-\r_n^j) \dt. \label{lwr_pi}
\eqn

In this section, we study the same problem as in Section 5.3, but with the kinematic wave model, \refe{ctm}, coupled with the discrete point queue model \refe{discretepq}. Here we set $\dx=30$ m, and the number of cells is $n=20$. Thus the CFL number $v_f \dt/\dx=1$. 
We consider both I- and PI-controllers in \refe{lwr_pi}. 

With the I-controller ($\beta=4$), the simulation results are given in \reff{vsl_pid_random_ctm}, where all solid lines are for the results with the I controller in \refe{lwr_pi}, and all dashed lines for those without VSL control; i.e., $u(t)=v_f$. From the figures, we can see that, when the arrival flow-rate is less than $C$ before $t=2000$ s, the queue size is 0, the in-flux equals the arrival flow-rate, as predicted by \refe{pointqueue}, and the traffic system has the same performance with or without control. However, at $t=2000$ s, the arrival flow-rate is greater than $C$ due to the disturbance, the I controller starts to reduce the speed limit in the VSL zone towards the optimal speed limit $v_1$, as shown in \reff{vsl_pid_random_ctm}(f), and the traffic density the last cell is controlled toward $k_1$, as shown in \reff{vsl_pid_random_ctm}(e). After that, the VSL control successfully increases the in-flux, as shown in \reff{vsl_pid_random_ctm}(c), prevents the occurrence of capacity drop and maintains a high discharging flow-rate, as shown in \reff{vsl_pid_random_ctm}(d), and maintains a small queue in the VSL and upstream region, as shown in \reff{vsl_pid_random_ctm}(b). 
Thus we can see that the PI controller is robust with respect to random variations in the demand patterns.
In this case, 2213 vehicles travel through the lane-drop region during 8000 s. The average travel time through the bottleneck is 39 seconds with the VSL control and 292 seconds without VSL control. Thus in this case the VSL control is able to save 86\% of the travel time. This confirms that the I-controller is also effective and robust for the kinematic wave model.

\bfg\bc
\includegraphics[width=5in]{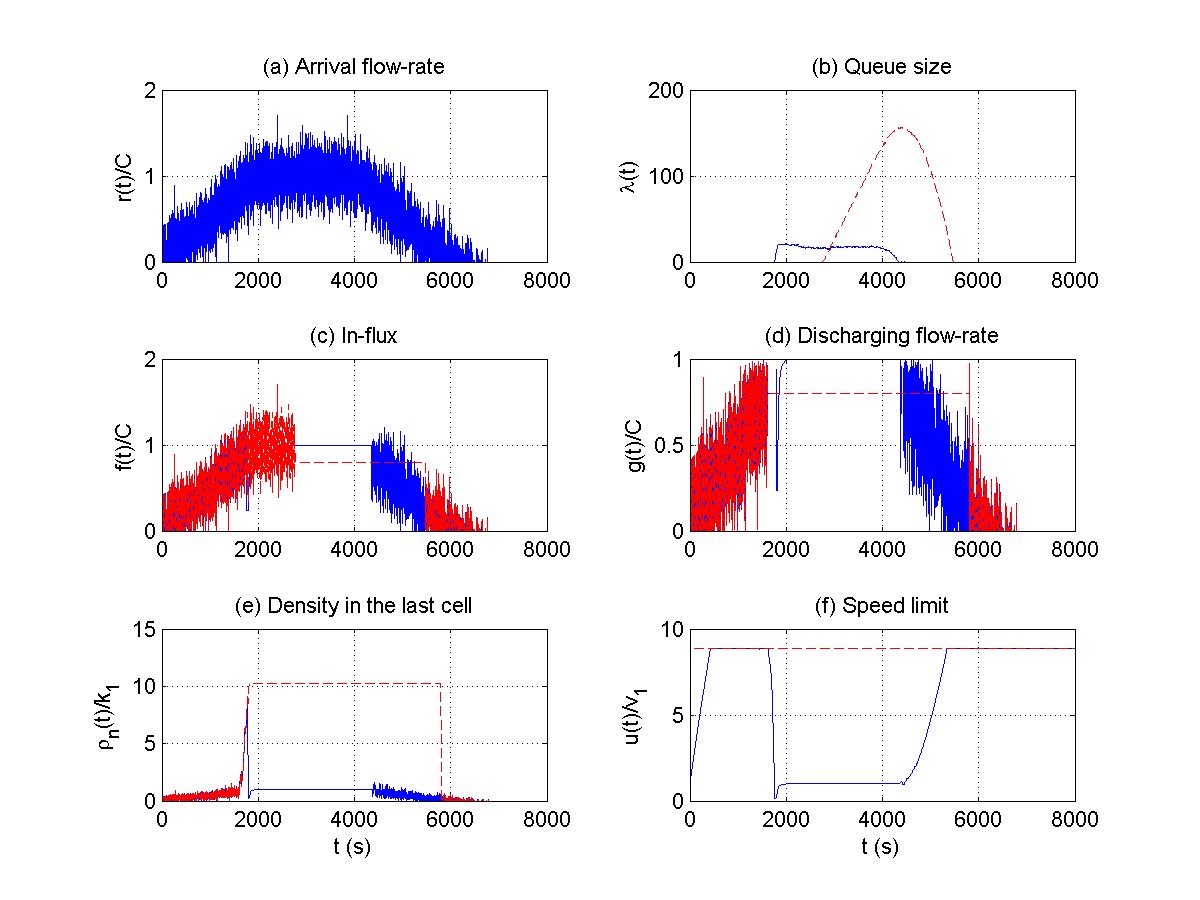}
\caption{Simulation of the kinematic wave model in a lane-drop region with (solid lines) and without (dashed lines) the VSL control}\label{vsl_pid_random_ctm}
\ec\efg

In \reff{vsl_pid_random_ctm_contour}, we compare the contour plot of traffic density $\r(x,t)$ inside zone 0 with and without the VSL control. From \reff{vsl_pid_random_ctm_contour}(b), we can see that, without the VSL control, the zone becomes congested at about 1800 s at the downstream boundary ($x=600$ m), the discharging flow-rate drops to $(1-\Delta)C$, and a shock wave forms and propagates upstream; the traffic queue starts to disappear at 5600 s along with a forward shock wave, when the upstream demand drops. In contrast, from \reff{vsl_pid_random_ctm_contour}(a), we can see that, with the VSL control, the zone also gets congested at 1800 s at $x=600$ m, but the I-controller immediately reduced the speed limit; then the system recovers to an uncongested state, and capacity drop is completely prevented.

\bfg\bc
\includegraphics[width=5in]{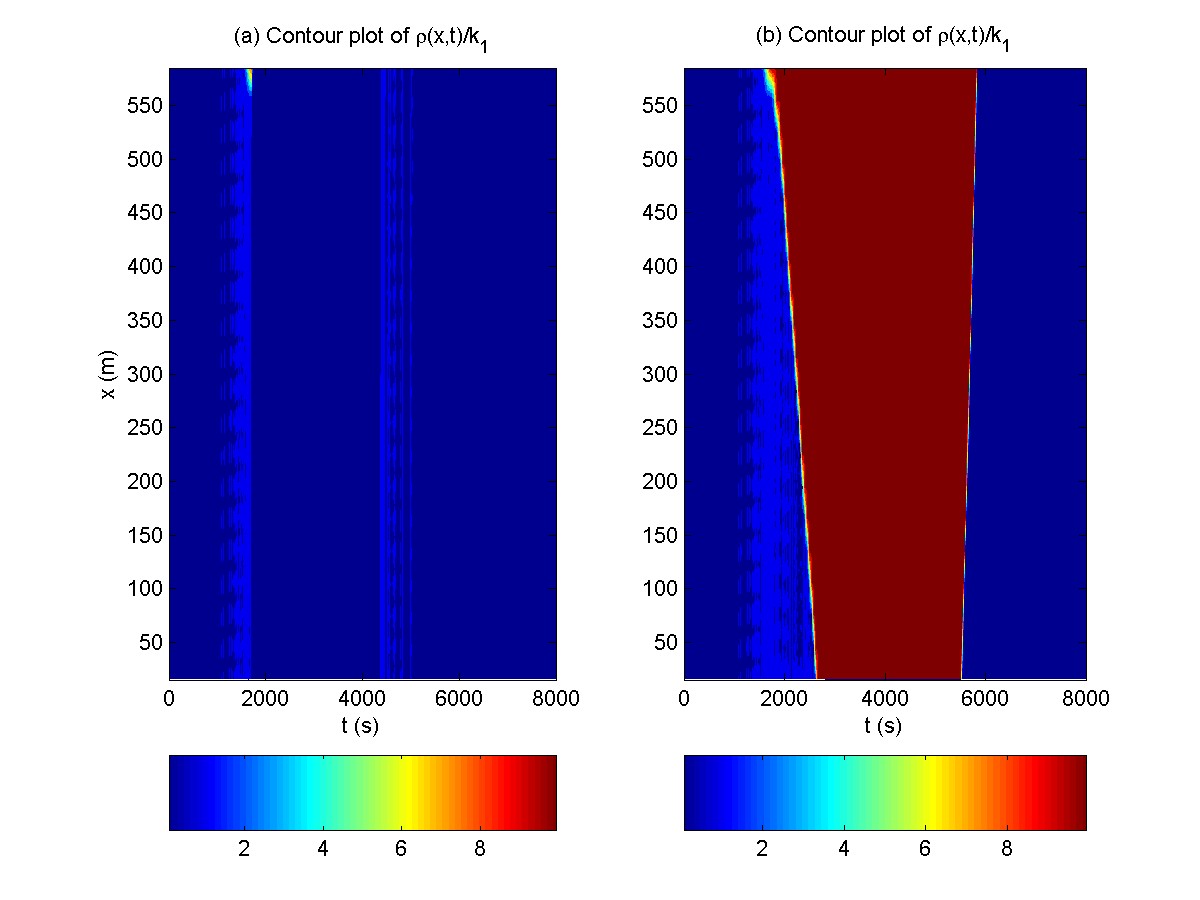}
\caption{Contour plot of traffic density for the kinematic wave model in a lane-drop region: (a) with the VSL control; (b) without control}\label{vsl_pid_random_ctm_contour}
\ec\efg

In \reff{vsl_pid_random_ctm_pi}, we demonstrate the simulation results with the PI controller ($\alpha=500$ and $\beta=20$). In this case, 2226 vehicles travel through the lane-drop region during 8000 s. The average travel time through the bottleneck is 43 seconds with the VSL control and 303 seconds without VSL control. Thus in this case the VSL control is able to save 86\% of the travel time. Therefore the PI-controller is also effective and robust for the kinematic wave model. Compared with results in \reff{vsl_pid_random_ctm}, the system performance is almost the same with the I- and PI-controllers.

\bfg\bc
\includegraphics[width=5in]{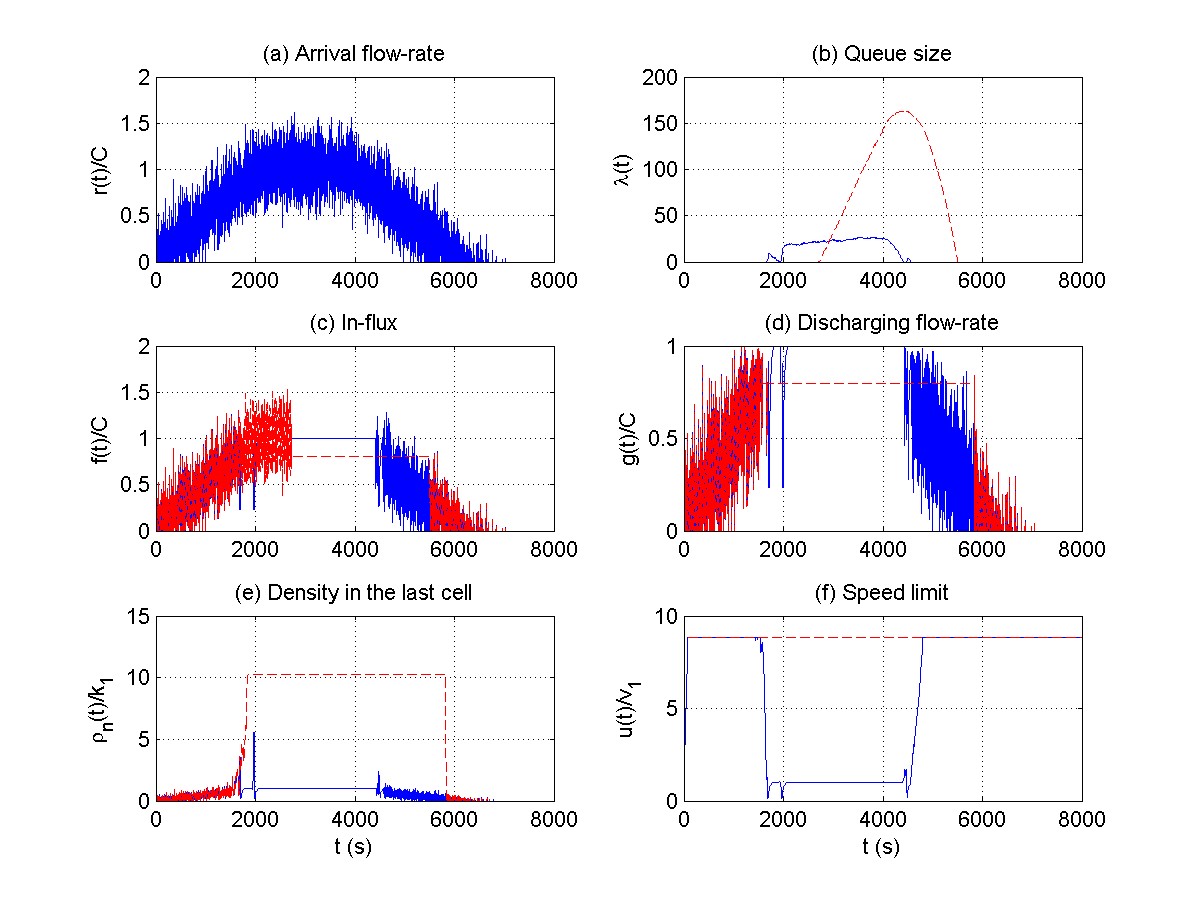}
\caption{Simulation of the kinematic wave model in a lane-drop region with (solid lines) and without (dashed lines) the VSL control: $\alpha=500$ and $\beta=20$}\label{vsl_pid_random_ctm_pi}
\ec\efg

\section{Conclusion}
In this paper, we studied the variable speed limits (VSL) control problem for a lane-drop bottleneck based on two traffic flow models: the kinematic wave model, a partial differential equation, and the link queue model, an ordinary differential equation. In both models, the discharging flow-rate is determined by a novel capacity drop model, and the in-flux is regulated by the speed limit in the VSL zone. Based on the link queue model, we first analytically proved that the open-loop control system with an optimal speed limit of $u^*=v_1$ has an uncongested equilibrium state with the maximum discharging flow-rate, but the other congested equilibrium state with capacity drop is always exponentially stable. Then we studied a state feedback control system with a target density $\bar k=k_1$ and proved that, with an I-controller, the congested equilibrium state is removed, and both I- and PI-controllers can stabilize the equilibrium state in the closed-loop control system. With numerical simulations, we demonstrated that the closed-loop system is stable with appropriate coefficients for the PI controller, the errors in the target density can substantially impact the performance of the system, and both I- and PI-controllers are robust subject to variances in the demand pattern. We also showed that 
With numerical simulations, we further showed that the I- and PI-controllers developed for the link queue model are also effective, stable, and robust for the LWR model.  

This study confirms that VSL strategies based on the PI- and I-controllers can effectively mitigate traffic congestion and reduce travel time when capacity drop occurs at a lane-drop bottleneck. On the other hand, analyses in Section 3.1 and simulations in Section 5.3 showed that, as expected, the VSL strategies cannot improve the performance of a traffic system without capacity drop, i.e., when $\Delta=0$. Therefore, before implementing VSL to improve traffic conditions, one needs to verify the existence of capacity drop at a study site.

Comparing results in Section 5.3 and Section 6, we find that the results are transferable between the link queue model and the LWR model. That is, a controller effective and robust for the link queue model can also be applied to control the LWR model. 
First, the performance of the control system is highly related to the stability of the switched dynamical systems in Section 4.2. In the future, we will be interested in studying their stability properties, quantifying the impacts of target densities, controller coefficients ($u_{min}$, $\alpha$, and $\beta$), traffic patterns, and other design parameters, including the length of zone 0, and determining whether discrete values of the speed limit can substantially change the equilibrium states and their stability.  
Second, theoretically this is a very interesting result, since the link queue model is a system of finite-dimensional ordinary differential equations, but the LWR model is a system of infinite-dimensional partial differential equations. Such a dual approach based on both the link queue and LWR models can be useful for studying other control problems of traffic flow, such as ramp metering and intersection signals. 
Third, the transferability of the VSL strategy from the link queue model to the kinematic wave model hints on the validity of the proposed strategy. In the future, we will be interested in examining whether the robust I-controller also works for more detailed car-following models and in the real world. 

In the future we will be interested in studying the impacts of various control strategies on the performance of the traffic system, such as bang-bang control and model predictive control.
It is also possible to adaptively change the target density under different road and traffic conditions. 
Another topic is to examine whether we can use the flow-rate instead of density as a state variable. In the LWR model in Section 6, we assume that the detector is installed in the last cell. A practically important topic is to study the impacts of the location of the traffic detector. We will also be interested in quantifying the environmental benefits of such control strategies for regular lane-drop bottlenecks, work zones, or incident areas \citep{smith2003characterization}.

\bibliographystyle{elsarticle-harv}

\end{document}